\documentclass[11pt]{article}

\usepackage[titletoc,title]{appendix}
\usepackage{amsfonts}
\usepackage{graphicx}
\usepackage{hyperref}
\usepackage{amsmath}
\usepackage{amssymb}
\usepackage{framed}
\usepackage{amsthm}
\usepackage{bbm}
\usepackage{tabstackengine}
\usepackage[usenames,dvipsnames,svgnames,table]{xcolor}

\newtheorem{theorem}{Theorem}[section]
\newtheorem{proposition}[theorem]{Proposition}
\newtheorem{lemma}[theorem]{Lemma}
\newtheorem{claim}[theorem]{Claim}

\newtheorem{definition}[theorem]{Definition}

\theoremstyle{remark}\newtheorem{remark}[theorem]{Remark}
\theoremstyle{remark}\newtheorem*{remark*}{Remark}

\newtheorem*{proposition*}{Proposition}
\theoremstyle{definition}\newtheorem{example}[theorem]{Example}

\usepackage{enumitem}
\newlist{casenv}{enumerate}{4}
\setlist[casenv]{leftmargin=*,align=left,widest={iiii}}
\setlist[casenv,1]{label={{\itshape\ \casename} \arabic*.},ref=\arabic*}
\setlist[casenv,2]{label={{\itshape\ \casename} \roman*.},ref=\roman*}
\setlist[casenv,3]{label={{\itshape\ \casename\ \alph*.}},ref=\alph*}
\setlist[casenv,4]{label={{\itshape\ \casename} \arabic*.},ref=\arabic*}
\providecommand{\casename}{Case}

\newcommand{\mrm}[1]{\mathrm{#1}}
\newcommand{\skipline}{$\phantom{}$}

\newcommand{\bem}[0]{\begin{eqnarray*}}
\newcommand{\enm}[0]{\end{eqnarray*}}
\newcommand{\anderbrace}[2]{
	\if\relax\detokenize{#2}\relax
		\sbox0{$\underbrace{#1}_{}$}
		\mathrel{\mathmakebox[\wd0]{#1}}
	\else
		\mathrel{\underbrace{#1}_{\mathclap{#2}}}
	\fi}

\newcommand{\integers}{\mathbb{Z}}
\newcommand{\pintegers}{\mathbb{N}}
\newcommand{\rationals}{\mathbb{Q}}
\newcommand{\reals}{\mathbb{R}}
\newcommand{\preals}{\reals_{\geq0}}
\newcommand{\spm}{\left\{-1,1\right\}}
\newcommand{\zo}{\left\{0,1\right\}}
\newcommand{\ozo}{\left\{-1,0,1\right\}}
\newcommand{\nicefrac}[2]{#1\left/#2\right.}

\newcommand{\sub}[0]{\subseteq}
\newcommand{\sm}[0]{\setminus}
\newcommand{\comp}{\circ}
\newcommand{\es}{\emptyset}
\newcommand{\ddd}{\doteqdot}
\newcommand{\im}{\mrm{Im}}
\newcommand{\one}{\mathbbm{1}}

\newcommand{\eps}{\epsilon}
\newcommand{\sgn}{\mrm{sgn}}

\newcommand{\pr}{\Pr}
\newcommand{\given}[2]{#1\,\middle|\, #2}
\newcommand{\normt}[1]{\lbn #1 \rbn_2}
\newcommand{\normts}[1]{\normt{#1}^{2}}
\newcommand{\vv}[1]{\mrm{Var}_{#1}}
\newcommand{\ee}[1]{\mathbb{E}_{#1}}
\DeclareMathOperator*{\be}{\mathbb{E}}
\DeclareMathOperator*{\var}{\mrm{Var}}

\newcommand{\andd}{\wedge}
\newcommand{\Andd}{\bigwedge}

\newcommand{\wh}[1]{\widehat{#1}}

\newcommand{\lbr}{\left(}
\newcommand{\rbr}{\right)}
\newcommand{\lbs}{\left[}
\newcommand{\rbs}{\right]}
\newcommand{\lbc}{\left\{}
\newcommand{\rbc}{\right\}}
\newcommand{\lba}{\left|}
\newcommand{\rba}{\right|}
\newcommand{\lbn}{\left\Vert}
\newcommand{\rbn}{\right\Vert}

\newcommand{\slice}[2]{\binom{[#1]}{#2 #1}}
\newcommand{\func}[3]{{#1}\colon {#2} \to {#3}}
\newcommand{\srfunc}[1]{\func{#1}{\slice{n}{p}}{\reals}}
\newcommand{\sbfunc}[1]{\func{#1}{\slice{n}{p}}{\zo}}
\newcommand{\sifunc}[1]{\func{#1}{\slice{n}{p}}{\integers}}
\newcommand{\AND}{\mrm{AND}}

\definecolor{ohadcolor}{RGB}{0, 100, 250}

\setstackEOL{\matEOL}
\newcommand{\restrict}[2]{{% we make the whole thing an ordinary symbol
  \left.\kern-\nulldelimiterspace % automatically resize the bar with \right
  #1 % the function
  \vphantom{\big|} % pretend it's a little taller at normal size
  \right|_{#2} % this is the delimiter
}}
\newcommand{\sspan}{\mrm{Span}}

\begin{document}

\title{A Structure Theorem for Almost Low-Degree Functions on the Slice}

\author{
	Nathan Keller\thanks{Department of Mathematics, Bar Ilan University, Ramat Gan, Israel.
		\texttt{nathan.keller27@gmail.com}. Research supported by the Israel Science Foundation (grants no.
		402/13 and 1612/17) and the Binational US-Israel Science Foundation (grant no. 2014290).}
	\mbox{ } and
	Ohad Klein\thanks{Department of Mathematics, Bar Ilan University, Ramat Gan, Israel.
		\texttt{ohadkel@gmail.com}.}
}

\maketitle

\begin{abstract}
The Fourier-Walsh expansion of a Boolean function $f \colon \zo^n \rightarrow \zo$ is its unique representation as a multilinear
polynomial. The Kindler-Safra theorem (2002) asserts that if in the expansion of $f$, the total weight on coefficients beyond degree
$k$ is very small, then $f$ can be approximated by a Boolean-valued function depending on at most $O(2^k)$ variables.

In this paper we prove a similar theorem for Boolean functions whose domain is the `slice' ${{[n]}\choose{pn}} = \{x \in \{0,1\}^n\colon \sum_i x_i = pn\}$,
where $0 \ll p \ll 1$, with respect to their unique representation as harmonic multilinear polynomials. We show that if in the representation of
$f\colon {{[n]}\choose{pn}} \rightarrow \zo$, the total weight beyond degree $k$ is at most $\epsilon$, where $\epsilon = \min(p, 1-p)^{O(k)}$, then $f$ can
be $O(\epsilon)$-approximated by a degree-$k$ Boolean function on the slice, which in turn depends on $O(2^{k})$ coordinates.
This proves a conjecture of Filmus, Kindler, Mossel, and Wimmer (2015). Our proof relies on hypercontractivity,
along with a novel kind of a shifting procedure.

In addition, we show that the approximation rate in the Kindler-Safra theorem can be improved from $\epsilon + \exp(O(k)) \eps^{1/4}$ to $\epsilon+\epsilon^2 (2\ln(1/\epsilon))^k/k!$, which is tight in terms of the dependence on $\eps$ and misses at most a factor of $2^{O(k)}$ in the lower-order term.
\end{abstract}

%%%%%%%%%%%%%%%%%%%%%%%%%%%%%%%%%%%%%%%%%%%%%%%%%%%%%%%
%%%%%%%%%%%%%%%%%%%%%%%%%%%%%%%%%%%%%%%%%%%%%%%%%%%%%%%
%%%%%%%%%%%%%%%%%%%%%%%%%%%%%%%%%%%%%%%%%%%%%%%%%%%%%%%
%%%%%%%%%%%%%%%%%%%%%%%%%%%%%%%%%%%%%%%%%%%%%%%%%%%%%%%
%%%%%%%%%%%%%%%%%%%%% INTRODUCTION %%%%%%%%%%%%%%%%%%%%
%%%%%%%%%%%%%%%%%%%%%%%%%%%%%%%%%%%%%%%%%%%%%%%%%%%%%%%
%%%%%%%%%%%%%%%%%%%%%%%%%%%%%%%%%%%%%%%%%%%%%%%%%%%%%%%
%%%%%%%%%%%%%%%%%%%%%%%%%%%%%%%%%%%%%%%%%%%%%%%%%%%%%%%
%%%%%%%%%%%%%%%%%%%%%%%%%%%%%%%%%%%%%%%%%%%%%%%%%%%%%%%

\section{Introduction}

\subsection{Background}

For a Boolean function $f:\{0,1\}^n \rightarrow \{0,1\}$, the Fourier-Walsh expansion of $f$ is its unique representation as an $n$-variate multilinear polynomial: $f(x) = \sum_{S \subset \{1,2,\ldots,n\}} \hat f(S) \chi_S(x)$, where $\chi_S(x)= \prod_{i \in S} (-1)^{x_i}$. The degree (or level) of a
coefficient $\hat f(S)$ is $|S|$, and a degree-$k$ function is a function for which $\hat f(S)$ vanishes for all $|S|>k$. The Fourier weight of
$f$ beyond degree $k$ is $W^{>k} (f) = \sum_{|S|>k} \hat f(S)^2$. (Note that by Parseval's identity, $\sum_{S} \hat f(S)^2$ is the expectation of $f^{2}$ with respect to the uniform measure on $\{0,1\}^n$.)

\paragraph{The Friedgut-Kalai-Naor (FKN) and the Kindler-Safra (KS) theorems.} The relations between the structure of a Boolean function and properties of its Fourier-Walsh expansion have been studied extensively in the last three decades. Many of the results achieved in this line of research rely on structural theorems characterizing Boolean functions whose Fourier-Walsh expansion has a `simple' form. The most basic of these is the Friedgut-Kalai-Naor (FKN) theorem~\cite{FKN}, which asserts that if most of the Fourier weight of $f$ lies on the two bottom levels, then $f$ essentially depends on a single coordinate.

We say that Boolean functions $f,g$ are $\epsilon$-close if $\Pr[f(x) \neq g(x)] \leq \epsilon$, where the probability is taken with respect to the uniform measure on $\{0,1\}^n$.
\begin{theorem}[Friedgut, Kalai, and Naor, 2002]
There exists a constant $C$ such that the following holds. Let $f$ be a Boolean function such that $W^{>1}(f) \leq \epsilon$. Then $f$ is $\epsilon$-close to one of the functions $0,1,x_i,1-x_i$ for some $1 \leq i \leq n$.
\end{theorem}
The FKN theorem has numerous extensions (see~\cite{ADFS,FF14,GH08,JOW15,MO10,Nayar14,O'Donnell14,Rub12}) and many applications, to hardness-of-approximation~\cite{Dinur07}, information theory~\cite{Samorodnitsky16}, social choice theory~\cite{FF14,Kalai-Choice}, extremal combinatorics~\cite{Friedgut08}, graph theory~\cite{ADFS,GH08}, and more.

\medskip
The Kindler-Safra (KS) theorem~\cite{KS} generalizes the FKN theorem to functions of a higher degree. It asserts that if most of the Fourier weight of $f$ is concentrated on the $k+1$ bottom levels, then $f$ can be approximated by a function that depends on at most $C^k$ coordinates.
\begin{theorem}[Kindler and Safra, 2002]\label{thm:original_ks}
There exist constants $C,c>0$ such that the following holds. Let $f$ be a Boolean function such that $W^{>k}(f) \leq \epsilon$, where $\epsilon < c^k$. Then $f$ is $(\epsilon+2^{Ck}\epsilon^{1/4})$-close to a Boolean-valued function depending on at most $C^{k}$ coordinates of $x$.
\end{theorem}
A recent argument of Dinur, Filmus and Harsha~\cite[Theorem 4.1]{DFH18} shows that for a sufficiently small constant $c$, the approximating function can be taken to be of degree at most $k$. In such a case, the approximating function depends on at most $6.614 \cdot 2^{k}$ coordinates, by a recent result of Chiarelli, Hatami, and Saks~\cite{CHS18} (which improved over a classical upper bound of $k \cdot 2^{k-1}$ by Nisan and Szegedy~\cite{NS94}).

We note that the Kindler-Safra theorem holds also for larger values of $\epsilon$, but the approximation rate becomes significantly weaker. In addition, it holds in the more general setting of the biased measure $\mu_p$ on the discrete cube; see~\cite{KS}.

Besides its intrinsic importance within the field of Boolean functions analysis, the Kindler-Safra theorem has a number of applications to extremal combinatorics (see~\cite{Friedgut08,Karpas17}).

\paragraph{Boolean functions on the slice.} While the initial results in Boolean functions analysis concerned properties of Boolean-valued functions on the discrete cube with respect to the uniform measure, multiple papers extended these results to functions over other domains, such as the solid cube $[0,1]^n$ (see~\cite{BKKKL}), the discrete cube endowed with a non-product measure (see~\cite{GG06}), and the symmetric group $\mathcal{S}_n$ (see~\cite{EFP11}).

One of these domains is \emph{a slice of the Boolean cube} (or `the slice' in short), composed of all elements of $\{0,1\}^n$ having the same Hamming weight.
\begin{definition}
Denote $[n]=\{1,2,\ldots,n\}$. For $0 \leq \ell \leq n$, the $\ell$'th slice of the discrete cube is $\binom{[n]}{\ell} \ddd \left\{ x \in \zo^n\colon \sum_{i=1}^n x_i = \ell \right\}$.
\end{definition}
Boolean functions over the slice are a natural object in hypergraph theory (where they correspond to properties of $\ell$-uniform hypergraphs on $n$ vertices), extremal combinatorics (where they correspond to properties of families of $\ell$-element subsets of a ground set of size $n$, which are the main object of study in intersection problems for finite sets, see~\cite{FT16}) and in coding theory (where they correspond to properties of constant-weight codes). In recent years, numerous papers studied Boolean functions over the slice (e.g.,~\cite{DT16,FKMW18,OW13,Wimmer14}). These papers generalized to the slice some of the most classical results of Boolean functions analysis and obtained applications to extremal combinatorics and to theoretical computer science.

In particular, O'Donnell and Wimmer~\cite{OW13} generalized the Kahn-Kalai-Linial (KKL) theorem~\cite{KKL}; Wimmer~\cite{Wimmer14} generalized Friedgut's junta theorem~\cite{Friedgut98} and Filmus~\cite{Filmus16} streamlined his proof using a new orthogonal basis for functions on the slice. Filmus~\cite{Filmus16b} generalized the FKN theorem~\cite{FKN}, and Das and Tran~\cite{DT16} used his result to settle a question of Bollob\'{a}s, Narayanan and Raigorodskii~\cite{BNR16} regarding the independence number of random subgraphs of the Kneser graph $K(n,k)$. Filmus et al.~\cite{FKMW18,FM16} generalized the Mossel-O'Donnell-Oleszkiewicz invariance principle~\cite{Maj-Stablest} and used it to obtain a generalization of a weak version of the Kindler-Safra theorem; this allowed them to obtain an improved stability version of the Ahlswede-Khachatrian theorem~\cite{AK} (in a certain range), which was later superseded in~\cite{EKL16} by different methods.

One of the main open problems raised by Filmus et al.~\cite[Problem~10.2]{FKMW18} was to prove a tight version of the Kindler-Safra theorem on the slice. In this paper we solve this problem.

\subsection{A generalization of the Kindler-Safra theorem to the slice}

In order to present our main result, we have to recall the generalization of the Fourier-Walsh expansion and of the notion of `degree' to functions on the slice.

It is clear that when the domain of the function is restricted to a slice of the discrete cube, i.e., when we consider $\func{f}{\binom{[n]}{\ell}}{\zo}$, the representation as a mulitilinear polynomial is not unique anymore. Instead, Srinivasan~\cite{Srinivasan11} and (independently) Filmus~\cite{Filmus16} suggested to use the representation as a \emph{harmonic} multilinear polynomial of degree $\leq \ell$ (assuming $\ell \leq n/2$). That is, the unique function
\[
g = \sum_{S\colon |S| \leq \ell} \hat g(S) \chi_S,
\]
where $g(x)=f(x)$ for all $x \in \binom{[n]}{\ell}$ and $\sum_{i \in [n]} \frac{\partial g}{\partial x_i} = 0$. Given the function $g$ that represents $f$, we denote
\[
	f^{>k} \ddd \sum_{S\colon |S| >k} \hat g(S) \chi_S
\]
and say that $f$ is of degree $\leq k$ if $f^{>k}=0$ (which means that $f$ can be represented by a degree-$k$ harmonic multilinear polynomial $g$).

In his generalization of the FKN theorem, Filmus~\cite{Filmus16b} shows that if a function $\func{f}{\slice{n}{p}}{\zo}$ satisfies
$||f^{>1}||_2^2 \leq \epsilon < cp^2$ for some universal constant $c$, then $f$ is $O(\epsilon)$-close to one of the functions $0,1,x_i,1-x_i$, for
some $1 \leq i \leq n$.\footnote{We note that the result of Filmus is more general, and applies also for larger values of $\epsilon$. In those cases (which are of interest mostly when $p$ is close to 0 or 1 --- a regime that we do not consider in this paper), one cannot obtain an approximation by such a simple function; Filmus shows that $f$ can be approximated by an affine function that depends on $O(\sqrt{\epsilon}/p)$ coordinates, which is the best one can obtain.}
Our main result is the following generalization of the Kindler-Safra theorem:\footnote{Throughout the paper, we assume for sake of simplicity that $pn$ is an integer. Of course, this does not affect the results.}
\begin{theorem}\label{thm:main}
There exists a constant $C$ such that the following holds. Let $0<p \leq 1/2$, and let $\sbfunc{f}$ satisfy $\normts{f^{>k}} \leq \eps \leq p^{Ck}$ (where the norm is taken with respect to the uniform measure on the slice). Then $f$ is $2\eps$-close to a degree-$k$ $\zo$-valued function $\tilde{f}$ on the slice.
\end{theorem}
Since, for $0 \ll p \ll 1$, degree-$k$ functions on the slice were shown in~\cite{FI18,FKMW18} to depend on $O(2^{k})$ coordinates, our theorem implies that $f$ essentially depends on $O(2^k)$ coordinates, similarly to the discrete cube case.

While Theorem~\ref{thm:main} does not imply the Kindler-Safra theorem directly, the argument can be modified
to yield a proof of the Kindler-Safra theorem. Alternatively, one can derive a slightly weaker version of the KS theorem by a `blackbox' reduction from Theorem~\ref{thm:main}, using the fact that the discrete cube $\{0,1\}^n$ with the biased measure $\mu_p$ can be `approximately embedded' into the slice $\slice{m}{p}$, for a sufficiently large $m$ (see, e.g.,~\cite[Theorem 3.3]{Filmus16b}). We present this reduction in Proposition~\ref{Prop:Reduction}.

\medskip We note that in~\cite[Theorem~8.5]{FKMW18}, Filmus et al. obtained a weaker version of Theorem~\ref{thm:main}, which makes the same hypothesis and shows that $\tilde{f}$ is $(\epsilon^{1/C'} + n^{-1/C'})$-close to $f$, for some constant $C'>0$. The authors of~\cite{FKMW18} conjectured that the approximation rate can be improved to $O(\eps)$, and this is proved in our Theorem~\ref{thm:main}.

\paragraph{Proof methods.} In general, our proof strategy is similar to the original proof of the Kindler-Safra theorem. In that proof, the approximating function $\tilde{f}$ is constructed via a sequence of $k+1$ functions, $\tilde{f}_k,\tilde{f}_{k-1},\ldots,\tilde{f}_0=\tilde{f}$, such that each function $\tilde{f}_l$ `well approximates' all discrete derivatives of $f$ of order $\geq l$. On the slice, discrete derivatives are replaced by operators $D_{ij}$ of the form $D_{ij}(f)(x)= (f(x)-f(x^{(ij)}))/2$, where $x^{(ij)}$ is obtained from $x$ by exchanging $x_i$ with $x_j$, and derivatives of higher order are obtained by sequential application of the derivative operator on pairwise disjoint pairs $(i,j)$ of coordinates.

We use a convenient-to-work-with basis of derivative operators, which we call `shifted sorted derivatives'. We show that each derivative operator can be represented as a linear combination of not-too-many shifted sorted derivatives, applied on permuted variants of the input. This allows us to construct a function which is guaranteed to `approximate well' only the shifted sorted derivatives and deduce that it approximates sufficiently well all other derivatives. The process of representing a derivative as a linear combination of shifted sorted derivatives resembles the classical combinatorial shifting technique (see~\cite{Frankl87}). We believe that this work strategy may be useful in other contexts as well.

In addition, the proof uses hypercontractive inequalities, mainly via the following lemma, whose `discrete cube' counterpart is the heart of the original argument of Kindler and Safra.
\begin{lemma}\label{lem:two_cases_norm}
There exists a constant $C$ such that the following holds. Let $0 <p \leq 1/2$ and let $\sifunc{f}$ and $k \in \pintegers$ be such that $\normts{f^{>k}}\leq \eps$. Then either $\normts{f}\leq 2\eps$ or $\normts{f} > p^{Ck}$.
\end{lemma}
The lemma asserts that any almost low-degree function from the slice to $\mathbb{Z}$ is either very close to the constant zero function, or attains non-zero values at a significant portion of its inputs. Together with other estimates, this lemma enables us to strengthen inequalities when they become loose, which turns out to be a crucial part of the proof strategy.

\paragraph{Future work.} One obvious approach for generalizing our results, is to allow a tradeoff between the strict guaranteed structure of the approximating function (being of degree $\leq k$), and the assumption $\epsilon < p^{O(k)}$, being too restrictive in some contexts. This approach was realized in previous works in different settings.
\begin{itemize}
\item Kindler and Safra~\cite{KS} proved a structure theorem similar to Theorem~\ref{thm:original_ks}, without the restriction on $\epsilon$; it is effective in the regime where $k=\omega(1)$. Though, expectedly, the guaranteed closeness to a junta is only $\epsilon^{1-\delta} \cdot c(\delta, p)$, with some $c(\delta, p) \to \infty$ as $\delta \to 0^{+}$.

\item Recently, Dinur, Filmus, and Harsha~\cite{DFH19} obtained a structure theorem for functions over the $p$-biased discrete cube, which is similar to Theorem~\ref{thm:original_ks}, but without the restriction on $\epsilon$ and with guaranteed closeness of $O(\epsilon)$; it is effective in the regime where $p=o(1)$. However, expectedly, the structure of the approximating function is greatly relaxed.
\end{itemize}
It is likely that both types of techniques can be combined with our methods, and yield a common generalization. We leave both approaches for future research.

\subsection{A sharpening of the Kindler-Safra theorem}

While the $O(\epsilon)$ approximation rate provided by the FKN and the Kindler-Safra theorems is usually sufficient, for some applications a more precise approximation rate is needed. For the FKN theorem, such a sharpening was obtained by Jendrej et al.~\cite{JOW15}, and independently by O'Donnell~\cite{O'Donnell14}:
\begin{proposition}[\cite{JOW15,O'Donnell14}]
There exists a constant $C$ such that the following holds. Let $f:\{0,1\}^n \rightarrow \{0,1\}$ be a Boolean function such that $W^{>1}(f) \leq \epsilon$. Then $f$ is $(\epsilon+C \epsilon^2 \log(1/\epsilon))$-close to one of the functions $0,1,x_i,1-x_i$ for some $1 \leq i \leq n$.
\end{proposition}
This sharpening was used by Samorodnitsky in his recent application of the FKN theorem to information theory~\cite{Samorodnitsky16}.

\medskip
We obtain a similar sharpening of the Kindler-Safra theorem:
\begin{theorem}\label{thm:sharp_ks}
There exists a constant $C > 0$ such that the following holds.
Let $\func{f}{\zo^{n}}{\zo}$ have $W^{>k} (f) \leq \eps \leq 2^{-Ck}$. Then, there exists a degree-$k$ function $\func{g}{\zo^{n}}{\zo}$ with
\begin{equation}\label{eq:sharp_ks_assertion}
\pr_{x}\lbs f(x) \neq g(x)\rbs \leq \eps + \eps^{2} \frac{(2\ln(1/\eps))^{k}}{k!}.
\end{equation}
\end{theorem}
We demonstrate, by an explicit example, that Theorem~\ref{thm:sharp_ks} is tight, up to a factor of $2^{O(k)}$ in the lower-order term.

\medskip We also show that a similar result holds for our generalization of the Kindler-Safra theorem to the slice.
\begin{theorem}\label{thm:sharp}
There exists a constant $C > 0$ such that the following holds. Let $0<p \leq 1/2$ and let $\sbfunc{f}$ satisfy $\normts{f^{>k}} \leq \eps \leq p^{C k}$ (where the norm is taken with respect to the uniform measure on the slice). Then $f$ is $\lbr \epsilon + \epsilon^2 \lbr C\log(1/p)\log(1/\epsilon)\rbr ^ {Ck} \rbr$-close to a degree-$k$ $\zo$-valued function $\tilde{f}$ on the slice.
\end{theorem}
% We also show that a similar sharpening holds in the `classical' discrete cube setting.
The proof is a bootstrapping over Theorem~\ref{thm:main}, using a generalization to the slice of the `level-$k$ inequalities' (see~\cite[Section~9.5]{O'Donnell14}).

\subsection{Organization of the paper}

The rest of the paper is organized as follows. In Section~\ref{sec:slice} we recall the basic notions of Boolean function analysis on the slice. In Section~\ref{sec:hyper-slice} we present the generalization of Bonami's hypercontractive inequality~\cite{Bonami} to the slice, and use it to prove Lemma~\ref{lem:two_cases_norm}. In Section~\ref{sec:derivatives} we introduce and study \emph{shifted sorted derivatives} that play a central role in the proof of the main theorem. The proof of Theorem~\ref{thm:main} is presented in Section~\ref{sec:main}, and we conclude with the proof of Theorems~\ref{thm:sharp_ks} and \ref{thm:sharp} in Section~\ref{sec:sharp}.

%%%%%%%%%%%%%%%%%%%%%%%%%%%%%%%%%%%%%%%%%%%%%%%%%%%%%%%
%%%%%%%%%%%%%%%%%%%%%%%%%%%%%%%%%%%%%%%%%%%%%%%%%%%%%%%
%%%%%%%%%%%%%%%%%%%%%%%%%%%%%%%%%%%%%%%%%%%%%%%%%%%%%%%
%%%%%%%%%%%%%%%%%%%%%%%%%%%%%%%%%%%%%%%%%%%%%%%%%%%%%%%
%%%%%%%%%%%%%%%%%%%%%%%% BASICS %%%%%%%%%%%%%%%%%%%%%%%
%%%%%%%%%%%%%%%%%%%%%%%%%%%%%%%%%%%%%%%%%%%%%%%%%%%%%%%
%%%%%%%%%%%%%%%%%%%%%%%%%%%%%%%%%%%%%%%%%%%%%%%%%%%%%%%
%%%%%%%%%%%%%%%%%%%%%%%%%%%%%%%%%%%%%%%%%%%%%%%%%%%%%%%
%%%%%%%%%%%%%%%%%%%%%%%%%%%%%%%%%%%%%%%%%%%%%%%%%%%%%%%

\section{Basics of Boolean Function Analysis on the Slice}
\label{sec:slice}

In this section we present some definitions and basic notions of analysis of Boolean functions on the slice that will be used in the sequel. For a more complete introduction, see~\cite{Filmus16,Filmus16b,FKMW18}.

\paragraph{Notation.} Throughout the paper, we use the following notations.
\begin{itemize}
\item
For $n \in \pintegers$, we let $[n]\ddd \{1,\ldots, n\}$.
%and $[-n]=\{-1,\ldots,-n\}$.
\item
For $I \subset [n]$, we denote the complement of $I$ by $\bar{I}=[n] \setminus I$.
\item
For $k \in \mathbb{N}$, we write $2^{-k} \mathbb{Z} \ddd \{2^{-k} m: m \in \mathbb{Z}\}$.
\item
For a family $S$ of sets, we denote $\bigcup S \ddd \lbc \given{x}{\exists T \in S:x\in T} \rbc$.
\item
For a finite set $A$, we write $a \sim A$ to mean that ``$a$ is a random variable uniformly distributed in $A$''.
\item
For a set $I$, we write $S_{I}$ for the set of permutations on $I$. For $n \in \pintegers$, we write $S_{n} \ddd S_{[n]}$.
\item
When we use a variable $x$, we always mean $x \in \zo^{n}$, and consequently, $x_{i}$ is a $\zo$-variable.
\item
For $n \in \pintegers$ and $p \in [0,1]$ such that $np \in \integers$, we write $\slice{n}{p}$ for the subset of $\zo^{n}$ specified by $\lbc \given{x \in \zo^{n}}{\sum x_{i} = pn} \rbc$. Usually, we will actually have $x \in \slice{n}{p}$.
\item
For a function $\srfunc{f}$ and for $q\in \mathbb{R}_{>0}$, we define the usual $L_{q}$ norm $\lbn f\rbn_{q}=(\be_{x}\lbs |f(x)|^{q} \rbs)^{1/q}$, where the expectation is taken with respect to the uniform measure on $\slice{n}{p}$.
\item
A random variable $X$ distributed according to the Poisson distribution with parameter $\lambda$ is described as $X \sim \mrm{Poi}(\lambda)$.

\item
By writing $a = O(b)$ (or $a=\Omega(b)$) we \emph{always} mean that there is a universal constant $c$, independent of any other parameter, such that $a \leq cb$ (or $b \leq ca$).
\end{itemize}

\paragraph{Assumption.}
We assume $0<p \leq \frac{1}{2}$, though it is implicit along the paper. All results can easily be extended to $p>1/2$ by replacing $p$ with $\min\{p,1-p\}$ throughout the paper. The only place this issue naturally emerges is in the introduction of `$\tau_{x}$' in the proof of Lemma~\ref{lem:eating_derivatives}.

\paragraph{Permutations, $k$-tuples, derivatives, and restrictions.} The following definitions will play a central role in the paper.
\begin{definition}[Permutations acting on the slice]
	Let $\pi \in S_{n}$ be a permutation. For an $x \in \slice{n}{p}$ we define $x^{\pi}$ as the element in $\slice{n}{p}$ satisfying
	\[
	\forall i \in [n]\colon (x^{\pi})_{\pi(i)} = x_{i},
	\]
or equivalently, $(x^{\pi})_{i} = x_{\pi^{-1}(i)}$.	For $\srfunc{f}$, we write $\srfunc{f^{\pi}}$ to describe the function $f^{\pi}(x)=f(x^{\pi})$.
\end{definition}
Note that the operator $f\mapsto f^{\pi}$ is linear, and that for $\pi, \tau \in S_{n}$, we have $(x^{\pi})^{\tau}=x^{\tau\comp\pi}$, and consequently, $(f^{\pi})^{\tau}=f^{\pi\comp\tau}$.

\begin{definition}[$k$-tuples]
	A set $P$ of $k$ disjoint ordered pairs from $[n]$ is called a \textbf{$k$-tuple}. More concretely, $P=\{(a_i, b_i)\}_{i \in [k]}$ is called a $k$-tuple, if
	\[
		\forall i\neq j \in [k]:\qquad
		a_{i} < b_{i},\qquad
		a_{i} \neq a_{j}, a_{i} \neq b_{j}, b_i \neq b_j, \qquad
		a_{i} \in [n], b_{i} \in [n].
	\]
	We say a $k$-tuple $P$ is \textbf{shifted} if
	\[
		\forall i:[a_{i}]\subseteq \bigcup P.
	\]
	$P$ is said to be \textbf{sorted} if
	\[
		a_{i}<a_{j} \Longleftrightarrow b_{i}<b_{j}.
	\]
We denote the set of all shifted sorted $k$-tuples by $\mathcal{V}_k = \mathcal{V}_k(n,p)$.
\end{definition}

\begin{definition}[Lexicographic order on the set of shifted sorted $k$-tuples]\label{Def:Lex}
For each $k \in \mathbb{N}$, we define a lexicographic order on $\mathcal{V}_k$, as follows. Let $P,P'$ be two shifted sorted $k$-tuples, and let $\{b_{i}\}, \{b'_{i}\}$ be the corresponding sorted sequences, satisfying $b_{i}\leq b_{i+1}$ and $b'_{i}\leq b'_{i+1}$ for all $i$. We say $P<P'$ if the minimal $j \in[k]$ with $b_{j} \neq b'_{j}$ has $b_{j} < b'_{j}$.
\end{definition}
Notice this is a total order, and no two distinct shifted sorted $k$-tuples are equivalent, as in a shifted-sorted $k$-tuple $\{(a_{i}, b_{i})\}_{i\in [k]}$, the $a_{i}$'s are uniquely determined by the $b_{i}$'s.

\begin{definition}[Derivative operator]\label{def:derivative_operator}
	For a $1$-tuple $P=\{(i,j)\}$, we define the corresponding \textbf{derivative operator} $D_{ij}$ by
	\[
	\forall f \in \reals^{\slice{n}{p}}\colon D_{ij}f \ddd \frac{1}{2}\lbr f - f^{(ij)} \rbr,
	\]
	where $(ij)$ is the transposition permutation $i\leftrightarrow j$. (Notice two disjoint $1$-tuples $p, q$ satisfy $D_{p}D_{q}=D_{q}D_{p}$.)	
	We generalize the \textbf{derivative operator} $D_{ij}$ to arbitrary $k$-tuples $P=\{(i_1,j_1),\ldots,(i_k,j_k)\}$, as
	\begin{equation}\label{eq:derivative_operator_definition}
	D_{P}f \ddd \be_{T\sub P} \lbs (-1)^{\lba T\rba} f^{T} \rbs = (D_{i_k j_k} \circ D_{i_{k-1} j_{k-1}} \circ \cdots \circ D_{i_1 j_1}) f.
	\end{equation}
	Furthermore, we say that a derivative $D_{P}$ is \textbf{shifted}, or \textbf{sorted}, if, correspondingly, $P$ is shifted, or sorted. We denote the set of all derivatives that correspond to shifted sorted $k$-tuples by $\mathcal{U}_k = \mathcal{U}_k(n,p)$.
\end{definition}
As happens in many similar contexts, it is not hard to see that $D_{ij}$ is a self-adjoint projection operator. We provide the simple proof of this property for the sake of completeness.
\begin{claim}\label{clm:projection_shrink}
	For any $1$-tuple $P=\{(i,j)\}$, $D_{ij}$ is an orthogonal projection. That is,
	\[
		\forall f,g\in \reals^{\slice{n}{p}}\colon \qquad D_{ij}(D_{ij} f) = D_{ij} f, \qquad and \qquad \left<D_{ij}f,g\right> = \left<f,D_{ij}g\right>.
	\]
	Consequently, $\normt{D_{ij} f} \leq \normt{f}$. Moreover, $D_{P}$ is an orthogonal projection for any $k$-tuple $P$, and consequently, $\normt{D_{P}f}\leq \normt{f}$ (i.e., $D_{P}$ is a contracting linear map).
\end{claim}
\begin{proof}
	First,
	\[
		D_{ij}(D_{ij} f) = \frac{(f - f^{(ij)})/2 - (f - f^{(ij)})^{(ij)}/2}{2} = (f-f^{(ij)}) / 2 = D_{ij} f.
	\]
	Then,
	\[
		2\left<D_{ij}f, g\right> = \left<f, g\right> - \be_{x}[f(x^{(ij)})g(x)] = \left<f, g\right> - \be_{y}[f(y)g(y^{(ij)})] = 2\left<f, D_{ij}g\right>,
	\]
	where the middle equality is due to $(ij)$ being an involution. Hence,
	\[
		\normts{f}-\normts{D_{ij}f}-\normts{f-D_{ij}f} = 2\left<D_{ij}f, f-D_{ij}f\right> = 2\left<f, D_{ij}f-D_{ij}(D_{ij}f) \right> = 0.
	\]
	Finally, $D_{P}$ is an orthogonal projection for any $k$-tuple $P$, from the commutativity of derivative operators corresponding to disjoint $1$-tuples.
\end{proof}	

\begin{remark}
We note that the inequality $\normt{D_{P}f}\leq \normt{f}$ can be easily derived using the triangle inequality. The stronger assertion that $D_P$ is an orthogonal projection will be used in the sequel.
\end{remark}

\begin{definition}[Restriction]
	For $x,y \in \slice{n}{p}$ and $I \subseteq [n]$, we write $x\equiv_{I} y$ to denote $\forall i\in I\colon x_{i}=y_{i}$.
	
\noindent
	Given a function $\srfunc{f}$, a set $I\sub [n]$ of coordinates, and $x \in \slice{n}{p}$, we write $f_{I=x}$ to describe the function which is the restriction of $f$ to the sub-slice $S=\lbc \given{y}{y \equiv_{I} x} \rbc \subseteq \slice{n}{p}$, i.e., $f_{I=x}=\restrict{f}{S}$. We further describe the function $\srfunc{\ee{I}f}$ by
	\[
	\ee{I}f(x) \ddd \be_{Q\in S_{I}} [f(x^Q)] = \be_{y}\lbs \given{f(y)}{y\equiv_{\bar{I}} x} \rbs = \be \lbs f_{\bar{I}=x} \rbs.
	\]
	Additionally,
	\[
		\vv{I}(f) \ddd \normts{f-\ee{I}f} = \be_{x} \lbs \var \lbr f_{\bar{I}=x} \rbr \rbs.
	\]
\end{definition}

\paragraph{Bases for the space of functions on the slice.}
In~\cite{Filmus16}, Filmus used the following definition towards a basis for the space of functions on the slice.
\begin{definition}
Let $Q$ be a $k$-tuple. The function $\func{\Psi_{Q}}{\slice{n}{p}}{\ozo}$ is defined by $\Psi_{Q}(x) = \prod_{(a,b)\in Q} (x_{b}-x_{a})$.
\end{definition}
Then, he explicitly built an appropriate basis.
\begin{definition}
Let $B \subseteq [n]$ have $|B|=k$. Define $\func{\chi_{B}}{\slice{n}{p}}{\integers}$ as the sum $\sum_{P} \Psi_{P}$, where $P$ ranges over all $k$-tuples $P=\{(a_i, b_i)\}_{i \in [k]}$ with $B = \{b_{1}, \ldots, b_{k}\}$.
\end{definition}
It turns out that the set $\{\chi_{B}\}$, where $B$ ranges over all subsets of $[n]$ of size $\leq pn$, for which $\chi_{B} \neq 0$ (that is, $B$ appears as $\{b_{1}, \ldots, b_{k}\}$ for any $k$-tuple), is a basis for the linear space of slice-functions $\reals^{\slice{n}{p}}$.

Instead of working with a specialized basis constructed for the space of slice-functions, we chose to use a basis for the derivative operators defined earlier. We believe this approach might be more natural in some contexts.
However, just in order to state the Kindler-Safra theorem on the slice, Theorem~\ref{thm:main}, one still needs a definition of the \emph{degree} of a function, and its \emph{decomposition} $f=\sum_{k} f^{=k}$ into the different levels. For this, we use a spanning set for the space of functions, which is common also in the context of the discrete cube $\{0,1\}^n$. In that context, the spanning set is $\lbc \AND_{S} \rbc_{S \subset [n]}$, where
\[
	\AND_{S} \ddd \Andd_{i \in S} x_{i} = \prod_{i \in S} x_{i}.
\]
Using this spanning set, one may define `the space of degree $\leq k$ functions on the discrete cube' to be $\sspan_{|S|\leq k} \lbc \AND_{S} \rbc$. This, of course, naturally expresses the set of multilinear polynomials of degree $\leq k$. This leads us to make a similar definition on the slice. Not coincidentally, it will turn out this definition coincides with that of Filmus~\cite{Filmus16}.
\begin{definition}\label{def:level_k}
For each $k \in \mathbb{N}$, `the space of degree $\leq k$ functions on the slice' is $L_{k} = \sspan_{|T|\leq k} \lbc \AND_{T} \rbc$, and `the space of homogeneous degree $k$ functions' is $R_{k} = L_{k} \cap L_{k-1}^{\perp}$. The `$k$-th level part' of a function $f$ on the slice, denoted $f^{=k}$, is the orthogonal projection of $f$ on $R_{k}$, with respect to the natural inner product on the space of functions defined over the slice $\left<f,g\right>\ddd\be_{x}[f(x)g(x)]$.
\end{definition}

The following claim shows that this alternative definition of $f^{=k}$ coincides with the definition of Filmus~\cite{Filmus16}. We note that we never explicitly use the expansion of $f$, and only implicitly use results regarding it, via Theorems~\ref{thm:noise_fourier} and~\ref{thm:hyper_slice} below, which were proved by Filmus with respect to an explicit basis. As these theorems, as well as all our assertions, depend only on the `levels' $\{f^{=k}\}_{k \in [pn]}$, the claim implies that we indeed can use our alternative basis instead of the basis of Filmus, without affecting the results.

The claim was proved by Filmus and Mossel in~\cite[Lemma 3.17]{FM16*} (see also~\cite[page~118]{FG89}). We provide a different proof for the sake of completeness.\footnote{Formally, in order to use our basis instead of the basis of Filmus, we have to show, in addition, that $\sspan_{Q:|Q|\leq k} \{\Psi_{Q}\} = \sspan_{B:|B|\leq k}\{\chi_B\}$. This follows from~\cite[Lemma~2.3 and Theorem~3.1]{Filmus16}.}
\begin{claim}\label{clm:and_vs_chi}
	For any $k \in \pintegers$,
	\[
		\sspan_{T:|T|\leq k} \{\AND_{T}\} = \sspan_{Q:|Q|\leq k} \{\Psi_{Q}\}.
	\]
\end{claim}
\begin{proof} We assume $k \leq np$, as both $\AND_{T}$ and $\Psi_{Q}$ are identically $0$ if $|Q| > np$ and $|T| > np$.

\noindent $(\supseteq)$ If $|Q|=k$, then $\Psi_{Q}$ is a linear combination of $2^{k}$ $\AND$ functions depending on $k$ variables.
	
\noindent $(\subseteq)$ If $|T|=k$, then $\AND_{T}$ is a polynomial of degree $k$ in $\sum_{i \in T} x_{i}$ on the slice, given by $\binom{\sum_{i \in T} x_{i}}{k}$. For the proof, we shall find a collection $\mathcal{Q}$ of $k$-tuples, so that the function $g=\sum_{Q\in \mathcal{Q}} \Psi_{Q}$ is also a polynomial in $\sum_{i \in T} x_{i}$ of degree exactly $k$ (and not less). This yields there is some $\alpha \in \rationals$ such that $h=\AND_{T}-\alpha \cdot g$ is a polynomial of degree $\leq k-1$ in $\sum_{i \in T} x_{i}$. In turn, this means $h$ is a linear combination of $\AND$ functions depending on $\leq k-1$ variables, which, by induction on $k$, are a linear combination of functions $\Psi_{Q}$ with $|Q|\leq k-1$. So overall this would give a representation of $\AND_{T}$ as a combinations of $\Psi_{Q}$ functions with $|Q| \leq k$, as required.

As proposed, we set $\mathcal{Q}$ to be the set of non-sorted $k$-tuples $Q=\{(a_{i}, b_{i})\}_{i\in [k]}$ satisfying $\{b_{i}\}_{i\in [k]} = T$. Notice that for this proof, we allow the (non-sorted) $k$-tuples $Q$ to have $a_{i} > b_{i}$. This can be fixed by swapping each pair having $a_{i} > b_{i}$ and negating the sign of $\Psi_{Q}$ accordingly.
Writing $X = \sum_{i\in T} x_{i}$, simple combinatorics give
\begin{equation}\label{eq:compatibility_aux}
	\sum_{Q\in \mathcal{Q}} \Psi_{Q}(x) = (-1)^{k-X}k! \frac{\binom{np-X}{k-X}\binom{n-np-k+X}{X}}{\binom{k}{X}} \ddd p(X).
\end{equation}
Since $X$ may attain all integer values between $0$ and $k$ and only these values, $p(t)$ can be uniquely interpolated as a degree-$k$ polynomial in $t$, $p(t)=\sum_{i=0}^{k} a_{i} t^{i}$. We wish to verify $a_{k} \neq 0$. For this, recall the following property of discrete differentiation:
\[
	k!a_{k} = \sum_{i=0}^{k} (-1)^{k-i} p(i) \binom{k}{i} \underbrace{=}_{\eqref{eq:compatibility_aux}} k! \sum_{i=0}^{k} \binom{np-i}{k-i} \binom{n-np-k+i}{i},
\]
evidently implying $a_{k} > 0$, as required.
\end{proof}

\paragraph{Simple properties of the derivative operators.} The following simple claim will be used several times in the sequel.
\begin{claim}\label{clm:derivative_operator_properties}
	Let $P$ be a $k$-tuple, $T \subseteq [n]$, and $\srfunc{f}$.
	\begin{enumerate}
		\item If $|T| < k$, then $D_{P}(\AND_{T}) \equiv 0$.
		
		\item If $|T|=k$, then $D_{P}(\AND_{T}) \equiv 0$, unless each pair of $P$ contains one element of $T$, in which case $D_{P}(\AND_{T}) \in \lbc \Psi_{P}/2^k, -\Psi_{P}/2^k \rbc$.
		
		\item $D_{P}(f)=D_{P}(f^{\geq k})$.
		
		\item If $m \in \pintegers$, then $D_{P}(f)^{> m} = D_{P}(f^{> m})^{> m}$.
	\end{enumerate}
\end{claim}
\begin{proof}\skipline
	\begin{enumerate}
	\item Consider a $p \in P$ for which $p\cap T = \es$. We may view $D_{P}(\AND_{T})$ as $D_{p}D_{P\sm \{p\}}(\AND_{T})$. Noting that $D_{P\sm \{p\}}(\AND_{T})$ does not, even syntactically, depend on the elements of $p$, we infer that exchanging the values of the coordinates of $p$, does not affect this function. Hence, $D_{p}D_{P\sm \{p\}}(\AND_{T}) = 0$.
	\item Unless each pair in $P$ contains exactly one element of $T$, we could find a $p \in P$ with $p \cap T = \es$ and proceed like the previous item. Moreover, one can verify that
	\[
		D_{\{i,j\}} (x_{i} \cdot g(x)) = \frac{x_{i}-x_{j}}{2} \cdot g(x)
	\]
	whenever $\func{g}{\slice{n}{p}}{\reals}$ is invariant to exchanging the values of $x_{i}, x_{j}$. Using this recursively, one can verify that if $T$ intersects every pair in $P$ in a single element
	\[
		D_{P}\AND_{T}=\pm \Psi_{P} \AND_{T \sm \bigcup P},
	\]
	where the $\pm$ sign depends on the parity of $\lba \lbc \given{p\in P}{p\cap T = \{\min(p)\}} \rbc \rba$.
	\item According to Definition~\ref{def:level_k}, we may write $f = \sum_{T:|T|<k} \alpha_{T} \AND_{T} + f^{\geq k}$. Recall that $D_{P}$ is linear, and so the first part of this claim ($D_{P}(\AND_{T})=0$) implies
	\[
		D_{P}(f) = D_{P}\lbr \sum_{T:|T|<k} \alpha_{T} \AND_{T} + f^{\geq k} \rbr = D_{P}(f^{\geq k}).
	\]
	\item In order to prove this item, it suffices to see that any $\srfunc{g}$ with $\deg(g) \leq m$ has $\deg(D_{P}(g)) \leq m$. This immediately follows from the linearity of $D_{P}$ and that $D_{P}(\AND_{T})$ is a linear combination of some $\AND_{S}$ functions with $|S|=|T|$.
	\end{enumerate}
\end{proof}

\section{Hypercontractivity on the Slice}
\label{sec:hyper-slice}

The classical noise operator $T_{\rho}$ defined over $\zo^{n}$ independently flips each coordinate of an input $x \in \zo^{n}$ with probability $(1-\rho)/2$. The noise operator, together with hypercontractivity results, have been successfully applied in many contexts, including in the proofs of the classical KKL theorem~\cite{KKL} and Friedgut's Junta theorem~\cite{Friedgut98}.

In order to understand the extension of the noise operator to the slice, we look on an equivalent definition of $T_{\rho}$. Instead of flipping each coordinate of $x$ with probability $(1-\rho)/2$, we may flip it back and forth $L$ times, where $L\sim \mrm{Poi}(\ln(1/\rho))$. It is also equivalent to just perform $t$ times a flip of a random coordinate of $x$, where $t\sim \mrm{Poi}(n \cdot \ln(1/\rho))$.

In the context of the slice, the natural replacement of flipping coordinates of an $x \in \slice{n}{p}$, is applying random transpositions $x \mapsto x^{(ij)}$ with random $i,j \in [n]$ ($i \neq j$), as this kind of operation leaves the distorted value of $x$ inside $\slice{n}{p}$. Specifically, the noise operator $H_{t}$ defined over the slice, applies random transpositions $x \mapsto x^{(ij)}$ to the input, a number of times which is distributed $\mrm{Poi}(t)$. Similarly to the $\zo^{n}$ context, $t=\Theta(n)$ corresponds to a `constant' noise rate.

\begin{definition}\label{def:noise_on_slice}
	Let $\srfunc{f}$ be a function defined on the slice. Define $L(f)=f-\nicefrac{\sum_{i<j} f^{(ij)}}{\binom{n}{2}}$, and
\[
H_t(f)=\exp(-tL) f= \sum_{l=0}^{\infty} \frac{t^{l}(-L)^l}{l!}f.
\]
Equivalently, $H_t(f)$ is the expectation of $f$ applied on the input after employing $\mrm{Poi}(t)$ random transpositions on the input.
\end{definition}
The first of the two following basic results is classical (see, e.g.,~\cite[Cor.~4.5 and Lem.~5.5]{Wimmer14}), and the other was proved by Lee and Yau~\cite[Thm.~5]{LY98}. The form in which they are written here is cited from~\cite{Filmus16}, Lemma~6.1 and Proposition~6.2, respectively.
\begin{theorem}\label{thm:noise_fourier}
		Let $t>0$ and $\alpha = \exp \lbr \frac{-2t}{n-1}\rbr$. Then $H_t(f) = \sum_{k=0}^{n} \alpha^{k-k(k-1)/n} f^{=k}$.
\end{theorem}
\begin{theorem}[\cite{LY98}] \label{thm:hyper_slice}
	The log-Sobolev constant $\rho$ related to the operator $L$ satisfies $\rho^{-1} = \Theta\lbr n \log{\frac{1}{p(1-p)}}\rbr$.
	Consequently, for any $1\leq \alpha \leq \beta \leq \infty$ with $\nicefrac{(\beta-1)}{(\alpha-1)}\leq \exp(2\rho t)$ and $\srfunc{f}$, we have $\lbn H_{t} f\rbn_{\beta} \leq \lbn f\rbn_{\alpha}$.
\end{theorem}

We shall use the following two hypercontractive inequalities. For the sake of completeness, we provide the simple proofs.
\begin{lemma}\label{lem:split_levels}
	Let $\srfunc{f}$, let $t\in \preals$, and let $k \in \pintegers$. Then
$$\normts{f} \leq \normts{f^{>k}}+ \exp \lbr \frac{4tk(n-k+1)}{n(n-1)} \rbr \normts{H_t f}.$$
\end{lemma}
\begin{proof}
	Let $\alpha = \exp \lbr \frac{-2t}{n-1} \rbr$ as in Theorem~\ref{thm:noise_fourier}, then
	\bem
	\normts{f}
	& =& \normts{f^{>k}}+\sum_{d=0}^{k} \normts{f^{=d}}
	\\
	& \leq &
	\normts{f^{>k}}+\alpha^{2(k(k-1)/n-k)}\sum_{d=0}^{k} \alpha^{2(d-d(d-1)/n)} \normts{f^{=d}}
	\\
	& \underbrace{\leq}_{\text{Thm.}~\ref{thm:noise_fourier}} &
	\normts{f^{>k}} + \alpha^{2(k(k-1)/n-k)} \normts{H_t f}.
	\enm
\end{proof}

The following proof is similar to the proof of~\cite[Theorem~9.21]{O'Donnell14}.
\begin{lemma}\label{lem:bonami_beckner_slice}
	Let $\srfunc{f}$ be of degree $\leq k$. Then $\be[f^8] \leq p^{-O(k)}\be[f^2]^4$.
\end{lemma}
\begin{proof}
Algebraically extend the definition of $H_{t}$ to support also $t\leq 0$, so that Theorem~\ref{thm:noise_fourier} still holds; that is,
\[
H_{t}f = \sum_{k=0}^{n} \exp(-2t(k-k(k-1)/n)/(n-1)) f^{=k}, \qquad \mbox{ for } t \leq 0.
\]
Let $t=\ln(7)/(2\rho)$. Since $H_{t}, H_{-t}$ both satisfy Theorem~\ref{thm:noise_fourier}, we have
	\begin{equation}\label{eq:bonami_almost_done}
		\lbn f \rbn_{8} = \lbn H_{t} H_{-t} f \rbn_{8} \underbrace{\leq}_{\text{Thm.}~\ref{thm:hyper_slice}} \lbn H_{-t} f \rbn_{2}.
	\end{equation}
By the definition of $H_{-t}f$, and since $f$ is of degree $\leq k$, we have $\lbn H_{-t} f\rbn_{2}^{2} \leq \exp\lbr \frac{2tk}{n-1}\rbr \normts{f}.$ Thus, raising Equation~\eqref{eq:bonami_almost_done} to the $8$-th power, we deduce
	\[
	\be[f^8] \leq \exp \lbr \frac{8tk}{n-1} \rbr \be[f^2]^{4} = p^{-O(k)} \be[f^2]^{4},
	\]
as asserted.
\end{proof}

Now we are ready to prove Lemma~\ref{lem:two_cases_norm}. Let us recall its statement.

\medskip \noindent \textbf{Lemma~\ref{lem:two_cases_norm}.}
There exists a constant $C$ such that the following holds. Let $0 <p \leq 1/2$ and let $\sifunc{f}$ and $k \in \pintegers$ be such that $\normts{f^{>k}}\leq \eps$. Then either $\normts{f}\leq 2\eps$ or $\normts{f} > p^{Ck}$.

\begin{proof}
	For $t=\frac{\ln(3)}{2\rho}$, where $\rho$ is as in Theorem~\ref{thm:hyper_slice}, we have
	\bem
	\normts{f}
	& \underbrace{\leq}_{\text{Lem.}~\ref{lem:split_levels}} &
	\normts{f^{>k}} + p^{-O(k)} \normts{H_{t}f}
	\\
	& \underbrace{\leq}_{\text{Thm.}~\ref{thm:hyper_slice}} &
	\eps + p^{-O(k)} \be[f^{4/3}]^{3/2}
	\\
	& \underbrace{\leq}_{f \in\integers} &
	\eps + p^{-O(k)} \be[f^2]^{3/2}.
	\enm
	Hence, $x=\be[f^2]$ satisfies $x-p^{-O(k)}x^{3/2}\leq\eps$, which implies that either $x\leq 2\eps$, or $x \geq p^{O(k)}$.
\end{proof}

Another tool which will come out handy in Section~\ref{sec:sharp} is the level-$k$ inequalities for slice functions. The proof is similar to a proof in~\cite[Section 9.5]{O'Donnell14}.
\begin{lemma}[Level-$k$ inequality]\label{lem:level_k_slice}
	There exists some universal constant $c_{0} > 0$ such that the following holds. Let $\func{f}{\slice{n}{p}}{\ozo}$ have $\be \lba f\rba = \epsilon$, then for any $k \in \pintegers$ with $\eps \leq p^{c_{0}k}$,
	\begin{equation}\label{eq:level_k_slice_stmt}
	\normts{f^{\leq k}} \leq \epsilon^{2}\lbr \frac{c_{0}\log(1/p)}{k} \log(1/\epsilon) \rbr^{c_{0}k/\log(1/p)}.
	\end{equation}
\end{lemma}
\begin{proof}
	Let $t>0$ to be determined, $\alpha = \exp\lbr \frac{-2t}{n-1} \rbr$, and $\gamma = 1+\exp(-2\rho t)$. Then
	\[
		\normts{f^{\leq k}}
		\underbrace{\leq}_{\mrm{Thm.}~\ref{thm:noise_fourier}}
		\alpha^{-2k} \normts{H_{t}(f)}
		\underbrace{\leq}_{\mrm{Thm.}~\ref{thm:hyper_slice}}
		\alpha^{-2k} \lbn f\rbn_{\gamma} ^ {2}
		\leq
		\alpha^{-2k}\epsilon^{2/\gamma}
		\leq
		\alpha^{-2k}\epsilon^{2(2-\gamma)}.
	\]
	Optimizing over $t$ to minimize $\alpha^{-2k} \epsilon^{4-2\gamma}$, one finds that the optimal $t$ satisfies
	\[
		\exp(-2\rho t)= \frac{k}{\rho(n-1)\ln(1/\epsilon)}.
	\]
(Note that this value of $t$ indeed satisfies $t \geq 0$, since we assumed $\epsilon \leq p^{c_{0}k}$; this is the only place in the proof where this assumption is used.) Let $r=\rho(n-1)$, which satisfies $r=\Theta(\log(1/p))$ due to Theorem~\ref{thm:hyper_slice}. We have
	\[
		\normts{f^{\leq k}} \leq \alpha^{-2k}\epsilon^{4-2\gamma} = \epsilon^{2} \lbr \frac{er \log(1/\eps)}{k}\rbr ^{2k/r},
	\]
	which implies the assertion of the theorem, provided $c_{0} \geq \max \lbc r / \log(1/p), 2\log(1/p)/r \rbc$ (which we may take, since $r=\Theta(\log(1/p))$, as noted above).
\end{proof}

\section{Shifted Sorted Derivatives}
\label{sec:derivatives}

In this section we prove two properties of shifted sorted derivatives that will play a central role in the proof of Theorem~\ref{thm:main}.

The first proposition asserts that each derivative of order $l$ can be written as a combination of a not-too-large number of shifted sorted derivatives of order $l$, applied on permuted variants of the input. It implies that in order to find a function $\tilde{f}$ that well approximates a given function $f$ with respect to all derivatives of order $l$, it is sufficient to find $\tilde{f}$ that approximates $f$ `very well' with respect to only the shifted sorted derivatives.

The second proposition complements the first one by showing how one can construct a low-degree function that approximates a given function with respect to all shifted sorted derivatives of some order $l$.

\begin{proposition}\label{lem:derivative_generating_set}
Let $\srfunc{f}$ and let $k \in \pintegers$. If $\normt{D_{P'}f}\leq \eps$ for every \textbf{shifted sorted} $k$-tuple $P'$, then:
	\begin{enumerate}
		\item For every $k$-tuple $P$ we have $\normt{D_{P}f}\leq 2^{10k^{2}}\eps$.
		\item If, furthermore, every derivative $D_{P}f$ satisfies that either $\normt{D_{P}f} \leq \eps$ or $\normt{D_{P}f} > 3\eps$, then all derivatives $D_{P}f$ satisfy $\normt{D_{P}f}\leq \eps$.
	\end{enumerate}
\end{proposition}

\begin{remark*}
We note that the bound $2^{10k^2}$ in the first assertion can be improved; however, we do not try to optimize it since in the rest of the paper we use only the second assertion.
\end{remark*}

In order to prove Proposition~\ref{lem:derivative_generating_set}, we perform a shifting procedure, in which a given derivative operator is gradually replaced by linear combinations of derivatives that are `closer to be shifted sorted'. To this end, we use the following two identities:
\begin{claim}
For any function $\func{f}{\slice{n}{p}}{\mathbb{R}}$ and for any $x\in \slice{n}{p}$, we have:
\begin{equation}\label{eq:derivative_flip}
(D_{12}D_{34}f)(x) = (D_{13} D_{24} f)(x^{(14)}) + (D_{14} D_{23} f)(x^{(24)}),
\end{equation}
and
\begin{equation}\label{eq:derivative_alter}
(D_{23}f)(x) = (D_{12}f)(x^{(13)})+(D_{13}f)(x^{(12)}).
\end{equation}	
\end{claim}
Identity~\eqref{eq:derivative_alter} is equivalent to
\begin{equation}\label{eq:derivative_alter2}
\sum_{\sigma \in S_{3}} \sgn(\sigma) f(x^{\sigma}) = 0.
\end{equation}
This holds for any $f$, since for any $x$ we have $x_{1} = x_{2}$ or $x_{2} = x_{3}$ or $x_{1} = x_{3}$, and so~\eqref{eq:derivative_alter2} contains only three of the terms $f(x^{\sigma})$, each appearing once with a `+' sign and once with a `-' sign. Identity~\eqref{eq:derivative_flip} seems to not have such a nice form, but can be easily verified by case analysis.
%Both identities can be easily verified by case analysis.

We note that in the proof of Proposition~\ref{lem:derivative_generating_set} below, identity~\eqref{eq:derivative_alter} can be replaced with the following identity:
\begin{equation}\label{eq:derivative_replacement}
(D_{23}f)(x) = (D_{12}f)(x)+(D_{13}f)(x^{(12)})+(D_{12}f)(x^{(123)}).
\end{equation}
Unlike~\eqref{eq:derivative_flip} and~\eqref{eq:derivative_alter}, Equation~\eqref{eq:derivative_replacement} does not assume $x\in \slice{n}{p}$, and may be applied whenever the ambient space (which is, in our case, $\slice{n}{p}$) is closed under the operation of permuting coordinates. In contrast, the identities arising from the formal expansions of~\eqref{eq:derivative_flip} and~\eqref{eq:derivative_alter} might look incorrect at first, unless one recalls $x_{i} \in \zo$.

\begin{proof}[Proof of Proposition~\ref{lem:derivative_generating_set}]
The proof uses the technique of `invariant'. We introduce a semi-invariant measure $m$ assigning a positive integer to every $k$-tuple, and show that given any non-shifted or non-sorted $k$-tuple $P$, and some permutation $\pi \in S_{n}$, we can express
\begin{equation}\label{eq:tuple_induction}
D_{P}f(x^\pi) = D_{S}f(x^\sigma) + D_{T}f(x^\tau),
\end{equation}
for some permutations $\sigma, \tau \in S_{n}$ and $k$-tuples $S,T$ such that $\max(m(S),m(T))< m(P)$. This, together with the triangle inequality, would inductively imply the second part of Proposition~\ref{lem:derivative_generating_set}. Furthermore, we will also show that our measure $m$ is upper bounded by $10k^2$, which would imply the first part of the Proposition. An example of the flow of the proof is presented after the proof.

\medskip We start by defining the semi-invariant measure $m$. For a $k$-tuple $Q$, consider its number of inversions
\[
\mrm{inv}(Q) = |\{((a,b),\,(a', b')) \in Q \times Q: (a < a') \wedge (b > b')\}|.
\]
We let $d(Q)$ be the sum of differences
\[
d(Q)=\sum_{(a,b) \in Q} \lbr i_{Q}(b)-i_{Q}(a) \rbr,
\]
where $i_{Q} \colon \bigcup Q \to [2k]$ is the rank-within-$Q$ map $i_{Q}(x) = \lba \lbc \given{y \in \bigcup Q}{y \leq x}  \rbc \rba$. In addition, we let
\[
I(Q) = (3k+1)\lbr 2k - \max \lbc \given{l \in [n]}{[l] \subseteq \bigcup Q} \rbc \rbr.
\]
Finally, $m(Q)$ is defined as
\[
m(Q) = \mrm{inv}(Q) + d(Q) + I(Q).
\]
We now show that given any non-shifted or non-sorted $k$-tuple $P$, and any permutation $\pi \in S_{n}$, we can express $D_P f(x^{\pi})$ in the form~\eqref{eq:tuple_induction}. We consider two cases:

\medskip \noindent \emph{Case 1: $P$ is non-sorted.} In this case, there exist $w<x<y<z$ so that $(w,z),\,(x,y) \in P$. Write $P' = P \sm \{(w,z),\,(x,y)\}$, so that an application of~\eqref{eq:derivative_flip} (with $w,z,x,y$ in place of $1,2,3,4$) yields:
\[
D_{P}f(x^\pi) = (D_{wz} D_{xy} D_{P'}f)(x^{\pi}) = (D_{wx} D_{zy}D_{P'}f)(x^{(wy)\pi}) + (D_{wy}D_{zx}D_{P'}f)(x^{(zy)\pi}).
\]
We have to verify that $S'=P'\cup \{(w, x),\,(y,z)\}$ and $T'=P'\cup \{(w, y),\,(x,z)\}$ satisfy $m(S') < m(P)$ and $m(T') < m(P)$. This is in turn implied from the following, which can be confirmed by case analysis:
\begin{enumerate}
\item $\mrm{inv}(T') < \mrm{inv}(P)$ and $d(T')=d(P)$.
\item $\mrm{inv}(S') + d(S') < \mrm{inv}(P) + d(P)$. (Note that $d$ is introduced as $\mrm{inv}(S')>\mrm{inv}(P)$ is possible.)
\item $I(T')=I(P)=I(S')$.
\end{enumerate}
Let us prove just the second assertion, which is the most tedious of the three. Clearly, one inversion appearing in $P$ but not in $S'$ is the one coming from $(w,z),\, (x,y)$. The inversions appearing in $S'$ and not in $P$ must have one pair in $P'$, and the other in $S'\sm P'$.

Let $(a,b)\in P'$. The number of inversions between $(a,b)$ and $P \sm P'$ (at most two inversions) is at least the number of inversions between $(a,b)$ and $S' \sm P'$, unless $\{a,b\}\cap \{x,y\} \neq \es$, in which case $(a,b)$ might be in at most one inversion in $S'$ (more than in $P$).

However, let us consider the difference $d(P)-d(S')$. Since $\bigcup S' = \bigcup P$, we have $i_{P} = i_{S'}$ and the sums of $d(P)$ and $d(S')$ differ only on the terms involving pairs from $S'\triangle P$. That is,
\[
d(P)-d(S')=2(i_{P}(z)-i_{P}(w)) \geq 0,
\]
and this quantity is at least twice the number of pairs $(a,b)\in P'$ with the mentioned property that increases the number of inversions of $S'$, by at most $1$. Overall, we observed at least one resolved inversion in $S'$ compared to $P$, and that every excess inversion in $S'$ corresponds to a reduction of at least $2$ in $d(S')$ compared to $d(P)$.

\medskip \noindent \emph{Case 2: $P$ is sorted but non-shifted.} In this case, there exists $(a,b) \in P$ so that $[a] \not \subseteq \bigcup P$. Let $a$ be minimal with these properties, and let $u=\min ([n]\sm \bigcup P)$. Write $P''=P \sm \{(a,b)\}$ so that an application of~\eqref{eq:derivative_alter} (with $u,a,b$ in place of $1,2,3$) yields:
\[
	D_{P}f(x^\pi) = (D_{ab} D_{P''}f)(x^{\pi}) = (D_{ua} D_{P''}f)(x^{(ub)\pi}) + (D_{ub} D_{P''}f)(x^{(ua)\pi}).
\]
We have to verify that $S''=P''\cup \{(u, a)\}$ and $T''=P''\cup \{(u, b)\}$ satisfy $m(S'') < m(P)$ and $m(T'') < m(P)$. For this we use the following claims. The first two claims state that a small variation in a $k$-tuple cannot increase too much its $d(\cdot)$ and $\mrm{inv}(\cdot)$ measures; the third states that our choice of $u$ to be the minimal unassigned number in $[n]$ decreases the $I(\cdot)$ measure by at least $3k+1$.
\begin{enumerate}
\item $d(S'') \leq d(P) + 2k$ and $\mrm{inv}(S'') \leq \mrm{inv}(P) + k$.
\item $d(T'') \leq d(P) + 2k$ and $\mrm{inv}(T'') \leq \mrm{inv}(P) + k$.
\item $[u] \subseteq \bigcup S''$, and since $u<a$, also $[u] \subseteq \bigcup T''$, while $[u] \not\subseteq \bigcup P$. Consequently, it holds that $\max\{I(S''), I(T'')\} \leq I(P) - (3k+1)$.
\end{enumerate}
\noindent All three claims are easy to verify.

\medskip Since $m(\cdot)$ of a $k$-tuple is non-negative, and since as long as $P$ is not shifted sorted, we can descend using~\eqref{eq:tuple_induction}, we can express each derivative as a linear combination of permuted shifted sorted derivatives. As for any $k$-tuple $P$, we have
\[
I(P) \leq 6k^2+2k, \qquad d(P) \leq k^2, \qquad \mbox{ and } \qquad \mrm{inv}(P) \leq \binom{k}{2},
\]
we overall have $m(P)\leq 9k^2$ for all $P$, which implies the first assertion of the Proposition.

The second assertion follows by induction on the measure of the $k$-tuple $P$. Indeed, using~\eqref{eq:tuple_induction} and the triangle inequality, $\normt{D_P f}$ can be bounded by the sum of the norms of two lower-measure derivatives of $f$. By the induction hypothesis, this means that $\normt{D_P f} < 2 \eps$. Together with the assumption that $\normt{D_P f} \not \in (\eps, 3\eps)$, we conclude $\normt{D_P f} < \eps$, as asserted.
\end{proof}

\begin{example}
For the convenience of the reader, we demonstrate the process described in the proof of Proposition~\ref{lem:derivative_generating_set}. Assume we start with the second-order derivative $D_{28} D_{67} f$ (i.e., $D_P f$, where $P=\{(2,8),(6,7)\}$). Let us compute $m(P)$. We have $\mrm{inv}(P) = 1$, $d(P) = 4$, and $I(P) = 7 \cdot (4-0) = 28$, and hence, $m(P) = 1 + 4 + 28 = 33$.

As $P$ is non-sorted, we may use~\eqref{eq:derivative_flip} to replace $D_P f$ by
\[
D_P f(x) = (D_{26}D_{78}f)(x^{(27)}) + (D_{27} D_{68} f)(x^{(78)}).
\]
A quick calculation shows that the measures $m$ of the new 2-tuples are
$m(\{(2,6),(7,8)\}) = 0 + 2 + 28 = 30$ and $m(\{(2,7),(6,8)\}) = 0 + 4 + 28 = 32$. Both are indeed smaller than $m(P)=33$.

Now, consider $Q=\{(2,6),(7,8)\}$, which is sorted. Using~\eqref{eq:derivative_alter}, we may write
\[
(D_{78} D_{26} f)(x) = (D_{17} D_{26} f)(x^{(18)}) + (D_{18} D_{26} f)(x^{(17)}).
\]
The measures $m$ of the two new 2-tuples are $m(\{(1,7),(2,6)\}) = 1 + 4 + 7 \cdot 2 = 19$ and $m(\{(1,8),(2,6)\}) = 1 + 4 + 7 \cdot 3 = 26$. Both are indeed smaller than $m(Q)=30$.

The process can be continued until all derivatives in the expansion of $D_P f$ as a sum, are sorted and shifted.
\end{example}

\begin{proposition}\label{lem:low_level_attain_derivatives}
For any $n,p,k$ and any function $z \colon \mathcal{V}_k \to 2^{-k}\integers$ (i.e., any function that assigns a value $\in 2^{-k}\integers$ for every shifted sorted $k$-tuple $P$ over the slice $\slice{n}{p}$), there is a function $\sifunc{f}$ of degree $\leq k$, having $D_{P} f = z(P) \Psi_{P}$ for every shifted sorted $P$.
\end{proposition}

Let us recall a relevant method underlying the proof we present. Suppose one is given some upper-triangular matrix $A\in \integers^{N\times N}$ satisfying $A_{ii}=1$ for all $i$. How could one efficiently express some $w\in \integers^{N}$ as an integral combination of the rows $a_{1},\ldots, a_{n}$ of $A$? The standard solution is to set $w_{0}=w$ and to iteratively define $w_{i}=w_{i-1}-(w_{i-1})_{i} a_{i}$ for $i=1,2,\ldots,N$, so that $w - \sum_{i=1}^{N} (w_{i-1})_{i} a_{i}= w_{N} = 0$.

That is analogous to how we shall proceed: First, we find an appropriate function-valued upper-triangular matrix whose rows correspond to all shifted sorted $k$-tuples, and whose columns correspond to shifted sorted derivatives, such that the diagonal entries are equal to the corresponding functions $\Psi_{P}$. Then, we construct the desired function $\sifunc{f}$ as an integral combination of the rows of the matrix. We demonstrate this process after the end of the proof.

\begin{proof}
For every shifted sorted $k$-tuple $P=\{(a_{i}, b_{i})\}_{i \in [k]}$ we let $B_{P}=\{b_{1},\ldots, b_{k}\}$ and define $g_{P} = \AND_{B_{P}}$.
	
Now, we arrange all the shifted sorted $k$-tuples according to the lexicographic order $<$ defined above (Definition~\ref{Def:Lex}) as $P_1 < P_2 < \cdots < P_m$, and consider the slice-function-valued matrix $C$ whose $i,j$ entry is $C_{ij}=D_{P_{j}}(g_{P_{i}})$.
	
The following three claims follow immediately from Claim~\ref{clm:derivative_operator_properties}(2).
	\begin{enumerate}
		\item $C$ is an upper-triangular matrix (i.e., below the diagonal we have zero functions).
        \item For each $i$, we have $C_{ii} = 2^{-k}\Psi_{P_{i}}$.
		\item Any non-zero entry $C_{ij}$ satisfies $C_{ij}=\pm 2^{-k}\Psi_{P_{j}}$.
	\end{enumerate}
Combining these claims, we see that if one takes $f_{0}=0$ and
\[
f_{i}=2^{k}\lbr z(P_{i})-t_{i} \rbr g_{P_i}+f_{i-1},
\]
where each $t_{i} \in 2^{-k}\integers$ is defined such that $D_{P_i}f_{i-1} = t_{i}\Psi_{P_i}$, then the final $f_{m}$ is the required function having the correct derivatives, i.e., $D_{P}f_{m} = z(P) \Psi_{P}$ for every shifted sorted $k$-tuple $P$.
\end{proof}

\begin{example}
Consider the case $n=5$ and $k=2$. In this case, there are five shifted-sorted 2-tuples:
\[
\{(1,2),(3,4)\}, \{(1,2),(3,5)\}, \{(1,3),(2,4)\}, \{(1,3),(2,5)\}, \mbox{ and } \{(1,4),(2,5)\}.
\]
The corresponding $\AND_{B}$ functions (sorted according to the lexicographic order $<$ defined above) are
\[
\AND_{(24)}, \AND_{(25)}, \AND_{(34)}, \AND_{(35)}, \mbox{ and } \AND_{(45)}.
\]
Let us construct the matrix $C$, as described in the proof of Proposition~\ref{lem:low_level_attain_derivatives}. As there are five possible $\AND_{B}$ functions, the matrix is of size $5 \times 5$. Its $(i,j)$ entry is $C_{ij}=D_{P_{j}}(g_{P_{i}})$, where the 2-tuples $P_i$ are arranged according to the lexicographic order in both rows and columns. A calculation shows that the resulting matrix is
\begin{equation*}
\setstacktabbedgap{-6pt}
\frac{1}{4}
\parenMatrixstack{
(x_2-x_1)(x_4-x_3)  & 0                     & 0                     & 0                     & (x_4-x_1)(x_2-x_5) \matEOL
0                     & (x_2-x_1)(x_5-x_3)  & 0                     & 0                     & 0 \matEOL
0                     & 0                     & (x_3-x_1)(x_4-x_2)  & 0                     & 0 \matEOL
0                     & 0                     &  0                    & (x_3-x_1)(x_5-x_2)  &   0 \matEOL
0                     & 0                     &  0                    &  0                    &   (x_4-x_1)(x_5-x_2)
}
.
\end{equation*}
Eventually, suppose we want to construct a function $f$ having desirable derivatives. We start by taking $\AND_{(24)}$ (i.e., the first row) to have the correct integer coefficient, so that the first derivative (i.e., $D_{\{(1,2),(3,4)\}}$, first column) of $f$ is the right multiple of $(x_2-x_1)(x_4-x_3)/4$. So, currently we have $f = c_1 \cdot \AND_{(24)}$, for some coefficient $c_1$. Then, we proceed to the second row $\AND_{(25)}$, and use the correct coefficient there to get the correct $D_{(1,2),(3,5)}$ derivative, and so on. Notice $C$ is not a diagonal matrix, so such an iterative procedure is required.
\end{example}

\section{Proof of Theorem~\ref{thm:main}}
\label{sec:main}

In this section we prove the main theorem of this paper. For the sake of convenience, we first present two lemmas which together imply the theorem (Lemmas~\ref{lem:eating_derivatives} and~\ref{lem:integer_to_boolean}) and prove the easier Lemma~\ref{lem:integer_to_boolean}. Then, in Section~\ref{sec:sub:lemma}, we present the more complex proof of Lemma~\ref{lem:eating_derivatives}.

Let us recall the formulation of Theorem~\ref{thm:main}.

\medskip \noindent \textbf{Theorem~\ref{thm:main}.} There exists a constant $C$ such that the following holds. Let $0<p \leq 1/2$, and let $\sbfunc{f}$ satisfy $\normts{f^{>k}} \leq \eps \leq p^{Ck}$. Then $f$ is $2\eps$-close to a degree-$k$ $\zo$-valued function $\tilde{f}$ on the slice.

\medskip
The main theorem easily follows from the combination of two lemmas.
\begin{lemma}\label{lem:eating_derivatives} There exists a constant $c_{1}$ such that the following holds.
Let $\sifunc{h}$ be a function with $\normts{h^{>k}}\leq \eps \leq p^{c_{1}k}$. Suppose there is some $l \leq k$ such that for every $(l+1)$-tuple $Q$ we have $\normts{D_{Q}h}\leq 2\eps$. Then there exists a function $\sifunc{g}$ of degree $\leq l$, such that for every $l$-tuple $P$ we have $\normts{D_{P}(h-g)}\leq 2\eps$.
\end{lemma}

\begin{lemma}\label{lem:integer_to_boolean}
	Let $\sifunc{f}$ be a function of degree $\leq k$. Either $f\equiv 0$ or $\pr[f \neq 0] \geq p^{O(k)}$.
\end{lemma}
We now deduce Theorem~\ref{thm:main} from the lemmas.
\begin{proof}[Proof of Theorem~\ref{thm:main}]
We shall recursively construct a sequence of functions $\sifunc{g_{k}, \ldots, g_{1}, g_{0}}$ such that for any $i$-tuple $P$ we have $\normts{D_{P}h_i} \leq 2\eps$, where $h_{i} = f-\sum_{j=i}^{k} g_{i}$.

For every $(k+1)$-tuple $P$ we know from Claim~\ref{clm:projection_shrink} and Claim~\ref{clm:derivative_operator_properties}(3) that
\[
	\normts{D_{P}h_{k+1}}=\normts{D_{P}(h_{k+1}^{>k})}\leq\normts{h_{k+1}^{>k}}\leq \eps,
\]
thus satisfying the hypothesis $\normts{D_{P}h_{k+1}} \leq 2\eps$. Lemma~\ref{lem:eating_derivatives} applied on $h_{i+1}$ recursively defines a $g_{i}$ which is of degree $\leq i$, so that for every derivative $P$ of order $i$ we have $\normts{D_{P} (h_{i+1}-g_{i})} \leq 2\eps$ (notice $h_{i}^{>k}=h_{i+1}^{>k}=f^{>k}$). Taking $g=\sum_{i=0}^{k} g_{i}$, we see that $\im(g)\subseteq \integers$ and
\[
	\normts{f-g} = \normts{D_{\es}(f-g)} \leq 2\eps.
\]
It only remains to show that $\im(g) \subseteq \zo$. This follows from Lemma~\ref{lem:integer_to_boolean}. Indeed, consider the function $r=g^2-g$ which is of degree $\leq 2k$. Notice that
\[
\pr[r \neq 0] = \pr[g \not \in \{0,1\}] \leq \pr[f\neq g]\leq 2\eps.
\]
Lemma~\ref{lem:integer_to_boolean}, applied on $r$, implies that either $r=0$ or $\pr[r\neq 0]\geq p^{O(2k)}$. However, we may assume $\eps$ is small enough so that the latter case cannot happen. Thus $g^2=g$, and so $\im(g) \subseteq \zo$.
\end{proof}
We note that Lemma~\ref{lem:integer_to_boolean} is a consequence of Lemma~\ref{lem:two_cases_norm} applied with $\epsilon = 0$. However, we present a tighter proof which is the `slice variant' of the one presented at \cite[Lemma~3.5]{O'Donnell14}.
\begin{proof}[Proof of Lemma~\ref{lem:integer_to_boolean}]
	We prove the lemma by induction on $n$ and $k$. Specifically, we prove that either $f \equiv 0$ or $\lba\{x\colon f(x)\neq 0\} \rba \geq \binom{n-2k}{np-k}$. (If $np < k$, we think of $\binom{n-2k}{np-k}$ as $1$.)
	
If $f$ is a constant function (which happens if $k=0$), either $f\equiv 0$ or $\pr[f\neq 0] = 1$, and the claim is true, so assume otherwise.
	
Consider $f'=D_{ij}f$ for some $i,j$ with $f'\not \equiv 0$. Such $i,j$ exist, as the graph defined over the slice $\slice{n}{p}$, where $(x\sim y) \iff \lbr \exists a,b\colon x^{(ab)}=y \rbr$, is connected. Since $g=f' / (x_{i}-x_{j})$ does not depend on $x_{i}, x_{j}$, it can be viewed as a nonzero degree $k-1$ polynomial defined on the $\binom{[n]\sm \{i,j\}}{np-1}$-slice
	\[
	\lbc\given{(x_{1},\ldots, \wh{x_{i}}, \ldots, \wh{x_{j}}, \ldots, x_{n})}{x\in\slice{n}{p} \andd x_{i} \neq x_{j}}\rbc.
	\]
Notice $g$ is indeed a degree $k-1$ polynomial, as $f'$ is a degree $k$ polynomial formally divisible by $x_{i}-x_{j}$. Moreover, $g$ does not depend on $x_{i}, x_{j}$ from the multilinearity of $f$.
	
By induction, $g(y)$ is nonzero for at least $\binom{n-2k}{np-k}$ values of $y$. For each such $y$, there is at least one $x$ with $f(x)\neq 0$, satisfying $\forall l\neq i,j\colon x_{l}=y_{l}$. Hence the induction is completed.

Overall, if $f \not \equiv 0$, we have
	\[
		\pr_{x} \lbs f(x)\neq 0\rbs \geq \nicefrac{\binom{n-2k}{np-k}}{\binom{n}{np}} \geq p^{O(k)}.
	\]
Notice this last inequality is true whenever $k \leq np/2$ by a simple computation using $(np-k)/n=\Omega(p)$, and for $k  > np/2$ since $\binom{n}{np} = p^{-O(k)}$.
\end{proof}

\subsection{Proof of Lemma~\ref{lem:eating_derivatives}}
\label{sec:sub:lemma}

In this subsection we prove Lemma~\ref{lem:eating_derivatives}, which is the core argument in the proof of Theorem~\ref{thm:main}, as was shown above.

To prove this lemma, we use three additional lemmas. The first is the following standard Chernoff-type bound for negatively correlated random variables, which follows from~\cite[Theorem~1]{PS97} via a classical Chernoff bound for $p$-biased Bernoulli variables.
\begin{lemma}\label{lem:neg_chernoff}
Let $x \in \slice{n}{p}$, and $S \subseteq [n]$. For any $t > 0$,
\[
	\pr \lbs \sum_{i \in S} x_{i} \leq p|S|-\sqrt{p(1-p)|S|}t \rbs \leq \exp(-t^{2}/2).
\]
\end{lemma}

The next lemma we use is a convenient form of Lemma~\ref{lem:two_cases_norm}. It provides a `dichotomy' statement regarding derivatives of integer-valued almost-degree-$k$ functions on the slice: each derivative of such a function is either close to the zero function, or attains a non-zero value `quite frequently'.
\begin{lemma}\label{lem:derivative_dichotomy}
There exists a real constant $c_{2}$ so that the following holds.
Let $\func{h}{\slice{n}{p}}{2^{-k}\integers}$ be a function with $\normts{h^{>k}} \leq \eps$. For any derivative $P$ of order $\leq k$ we have either $\normts{D_{P}h} \leq 2\eps$, or else $\normts{D_{P}h} \geq p^{c_{2}k}$.
\end{lemma}
\begin{proof}
	Write $H=D_{P} h$, for a derivative $P$ of order $\leq k$.
	Note that
	\begin{equation}\label{Eq:Aux5.00}
		\normts{H^{>k}}\leq\normts{D_{P}(h^{>k})}\leq \normts{h^{>k}}.
	\end{equation}
	The first inequality follows from Claim~\ref{clm:derivative_operator_properties}(4), that is, $D_{P}$ is linear and cannot increase the degree of its input, and so $H^{>k}=D_{P}(h^{>k})^{>k}$. The second inequality follows from Claim~\ref{clm:projection_shrink}, stating that $D_{P}$ is a contracting projection.
	
	Since $H$ is a derivative of the $2^{-k}\integers$-valued function $h$, we deduce that $H$ is $2^{-k-|P|}\integers$-valued. In particular, applying Lemma~\ref{lem:two_cases_norm} to the integral-valued function $2^{2k}H$ we find that either $\normts{2^{2k}H}\leq 2\cdot 2^{4k}\eps$ or $\normts{2^{2k}H} \geq p^{O(k)}$. The former case means $\normts{H}\leq 2\eps$ and the latter case means $\normts{H} \geq p^{c_{2}k}$, for some universal constant $c_{2}$.
\end{proof}
The third lemma we use concerns the conditional variance $\vv{I}(f)$. Recall that given $I \subseteq [n]$, $\vv{I}(f)$ measures the expected uncertainty in the value of $f(x)$ given only $x_{\bar{I}}$, for a random $x \in \slice{n}{p}$. We thus expect that if the value of a coordinate $x_{i}$ does not greatly influence $f(x)$, then $\vv{I}(f)$ should not be much larger than $\vv{I \setminus \{i\}}(f)$. The lemma confirms this intuition.
\begin{lemma}\label{lem:var_bound_by_larger}
	Let $\srfunc{f}$ be a function. Let $I \subseteq [n]$ be a set of coordinates, and let $i\in I$ satisfy that for any $j \in I$ we have $\normts{D_{ij}f} \leq \eps$. Then
	\[
		\vv{I}(f) \leq \vv{I\sm \{i\}}(f) + \eps.
	\]
\end{lemma}
\begin{proof}
Write $J=I \sm \{i\}$. For a uniform $x\sim \slice{n}{p}$, consider the two-step martingale
\[
	\ee{I}f(x), \ee{J}f(x), f(x).
\]
Since the expected value of $f(x)$, given $\ee{I}f(x)$ and $\ee{J}f(x)$, is $\ee{J}f(x)$, %(a martingale),
we have
\[
	\vv{I}f = \normts{f-\ee{I}f} = \normts{f-\ee{J}f} + \normts{\ee{J}f-\ee{I}f} = \vv{J}(f) + \normts{\ee{J}f-\ee{I}f}.
\]
Hence, it only remains to prove
\begin{equation}\label{Eq:Aux5.1}
\normts{\ee{J}f-\ee{I}f} \leq \eps.
\end{equation}
To this end, we claim that
\[
	\normts{\ee{J}f - \ee{I}f} \leq \be_{j \sim I} \lbs \normts{D_{ij}f} \rbs,
\]
where $D_{ii}=0$. To see this, define the operator $\func{T}{\reals^{\slice{n}{p}}}{\reals^{\slice{n}{p}}}$ by
\[
T(f) (x)=\be_{j \sim I} \lbs D_{ij}f(x) \rbs = f(x)-\be_{j \sim I}\lbs f^{(ij)}(x)\rbs.
\]
It is standard
that if $j\sim I$ and $Q\sim S_{I\sm \{i\}}$, then $(i, j)\comp Q$ is a uniform permutation in $S_{I}$. Thus,
\[
	\ee{J}(T(f)) = \ee{J}(f) - \ee{I}(f).
\]
Notice that $\ee{J}$ is an averaging operator, and additionally $T$ is an average of some operators $D_{ij}$. So, applying the inequality $\be[Z]^2\leq \be[Z^2]$ twice on these averagings implies
\begin{align*}
\normts{\ee{J}f - \ee{I}f} &= \be_{x}[\lbr \ee{J}f(x) - \ee{I}f(x) \rbr^{2}] \leq
	\be_{x}\lbs\ee{J}\lbs (Tf)(x)^{2} \rbs\rbs \\
&\leq \be_{x} \lbs \be_{j \sim I} \lbs (D_{ij} f(x))^{2} \rbs \rbs =
	\be_{j \sim I} \lbs \normts{D_{ij}f} \rbs.
\end{align*}
However, by assumption, for every $j\in I$ we have $\normts{D_{ij}f} \leq \eps$, and hence,~\eqref{Eq:Aux5.1} follows.
\end{proof}

Now we are ready to present the proof of Lemma~\ref{lem:eating_derivatives}.

\begin{proof}[Proof of Lemma~\ref{lem:eating_derivatives}]
The proof is composed of two steps. In the first step, we use hypercontractive estimates, together with Lemma~\ref{lem:var_bound_by_larger}, to show that for every derivative $P$ of order $l$, we have $\vv{J}(D_{P}h) \leq O(k^2/p)\eps$, where $J=[n]\sm \bigcup P$. This means that $D_{P}h$ can be approximated by a function of the form $a\cdot \Psi_{P}$ for some $a \in \reals$. In the second step we build the required function $g$. We choose $g$ so that all its shifted sorted derivatives of order $l$ approximate the corresponding ones of $h$; we do this via Proposition~\ref{lem:low_level_attain_derivatives}. Then, we use the second part of Proposition~\ref{lem:derivative_generating_set} to show that actually $g$ approximately agrees with $h$ on all order-$l$ derivatives, thus concluding the proof.

\medskip \noindent \textbf{Step~1. A bound on $\vv{J}(D_{P}h)$.} Let $P$ be a derivative of order $l$, and denote $H=D_{P}(h)$. Let $I\sub [n]$ satisfy $I \cap \bigcup P = \es$.
We shall prove by induction on $|I|$ that
\begin{equation}\label{eq:eating_var_bound}
	\vv{I}(H) \leq \min(K, 2|I|)\eps,
\end{equation}
where $K$ is a constant depending on $k,p$, to be determined later.

Denote $X = \lbc \given{x\in \slice{n}{p}}{\var \lbr H_{\bar{I}=x}\rbr > \tau_x}\rbc$, where $\tau_{x} = \nu \lbr \sum_{i \in I} x_{i}/|I|\rbr^{c_{2}k} / 5$, $\nu(t)=\min\{t,1-t\}$, and $c_{2}$ is the universal constant from Lemma~\ref{lem:derivative_dichotomy}. Notice that $\tau_{x}$ depends only on $x_{\bar{I}}$, and thus, $X$ is a union of sub-slices of $\slice{n}{p}$.

\medskip We clearly have
\begin{equation}\label{eq:eating_bound0}
	\mathrm{Var}_I(H) = \be_{x}\left[(H(x)-\ee{I}H(x))^{2}\one_{x\notin X}\right] + \be_{x}\left[(H(x)-\ee{I}H(x))^{2}\one_{x \in X}\right].
\end{equation}
We bound each term of the RHS separately. Unfortunately, this bounding argument --- which spans the following three pages --- is quite cumbersome.

\medskip \noindent \emph{Sub-step 1: Bounding $\be_{x}\left[(H(x)-\ee{I}H(x))^{2}\one_{x\notin X}\right]$.}
In order to bound this term, we replace $\ee{I}H(x)$ by an integral multiple of $2^{-k}$, depending only on $x_{\bar{I}}$, without increasing the expression we want to bound much. We let
\[
G(x)=2^{-k}\lfloor 2^{k} \ee{I}H(x)\rceil,
\]
where $t \mapsto \lfloor t \rceil$ is the rounding-to-nearest-integer function. Then on the one hand, we have
\begin{equation}\label{Eq:Aux5.2}
	\be_{x} \lbs (H(x)-\ee{I}H(x))^{2} \one_{x\notin X} \rbs \leq \be_{x} \lbs (H(x)-G(x))^2 \one_{x\notin X} \rbs,
\end{equation}
because of the well-known property that for every real-valued variable $A$, the function $a\mapsto\be[(A-a)^2]$ is minimized at $a=\be[A]$. On the other hand, we have
\begin{equation}\label{eq:round_not_so_bad}
	\forall x \in \slice{n}{p}\colon \lba H(x) - G(x) \rba \leq 2\lba H(x) - \ee{I}H(x) \rba,
\end{equation}
as $H$ is an order-$l$ derivative of $h$, and hence $\im(H) \subseteq 2^{-l}\integers \subseteq 2^{-k}\integers$.

\medskip Fix some $x\notin X$. Applying Lemma~\ref{lem:derivative_dichotomy} to the function $F=(H-G)_{\bar{I}=x}$, we infer that either $\normts{F} \leq 2\normts{F^{>k}}$ or $\normts{F} > 4\tau_{x}$ (recall the definition of $\tau_{x}$).
Notice we know from~\eqref{eq:round_not_so_bad} that
\[
\normts{F}\leq 4\normts{{(H-\ee{I}H)}_{\bar{I}=x}}=4\var \lbr H_{\bar{I}=x}\rbr.
\]
Since $x \notin X$, this implies $\normts{F}\leq 4\tau_{x}$. Henceforth, we must have
\begin{equation}\label{Eq:Aux5.3}
\normts{F} \leq 2\normts{F^{>k}}.
\end{equation}
In order to proceed, we claim
\begin{equation}\label{eq:eating_shrink}
	\normts{F^{>k}} \leq \normts{(H^{>k})_{\bar{I}=x}}.
\end{equation}
To see this, we observe that
\begin{equation}\label{Aux5.1}
F^{>k}=\lbr H_{\bar{I}=x} \rbr^{>k}=\lbr\lbr H^{>k} \rbr_{\bar{I}=x}\rbr^{>k}.
\end{equation}
The left equality is a combination of $F$'s definition and the fact that $G_{\bar{I}=x}$ is constant (as $y\mapsto G(y)$ depends only on $y_{\bar{I}}$). The right equality follows from the fact that restricting a function to a sub-slice does not increase its degree. Indeed, any $\AND_{S}$ defined on the slice $\slice{n}{p}$, when viewed on the sub-slice $\bar{I}=x$, remains of degree $\leq |S|$. Since $D\mapsto D^{>k}$ is a contracting projection (Claim~\ref{clm:projection_shrink}),~\eqref{Aux5.1} implies~\eqref{eq:eating_shrink}. As Equations~\eqref{Eq:Aux5.3} and~\eqref{eq:eating_shrink} hold for any $x \notin X$, it follows by averaging over all $x \notin X$ that
\begin{equation}\label{Eq:Aux5.4}
\be_{x} \lbs (H(x)-G(x))^2 \one_{x\notin X} \rbs \leq 2\normts{H^{>k}(x)\cdot\one_{x\notin X}}.
\end{equation}
Combining Equations~\eqref{Eq:Aux5.2} and~\eqref{Eq:Aux5.4}, we obtain
\begin{equation}\label{Eq:Aux5.5}
	\be_{x} \lbs (H(x)-\ee{I}H(x))^{2} \one_{x\notin X} \rbs \leq \be_{x} \lbs (H(x)-G(x))^2 \one_{x\notin X} \rbs \leq 2\normts{H^{>k}(x)\cdot\one_{x\notin X}} \leq 2\eps,
\end{equation}
where the last inequality follows from the assumption $||h^{>k}||_2^2 \leq \eps$ and~\eqref{Eq:Aux5.00}. This completes the first sub-step.

\medskip \noindent \emph{Sub-step 2: Bounding $\be_{x}\left[(H(x)-\ee{I}H(x))^{2}\one_{x\in X}\right]$.}
Since for any $a,b \in \mathbb{R}$, we have $(a+b)^2 \leq 2(a^2+b^2)$, it follows that
\begin{align}\label{eq:eating3}
\begin{split}
	\be_{x} \lbs \lbr H(x) - \ee{I}H(x) \rbr ^{2} \one_{x\in X}\rbs
	& =
	\be_{x} \lbs \lbr H(x) - \ee{I}H(x) \rbr ^{2} \one_{x\in X \andd \tau_{x} \neq 0}\rbs\\
	&\leq2 \be_{x} \lbs \lbr H^{\leq k}(x) - (\ee{I}H)^{\leq k}(x) \rbr ^ {2} \one_{x\in X \andd \tau_{x} \neq 0}\rbs
	\\
	&\ \ +2 \be_{x} \lbs \lbr H^{> k}(x) - (\ee{I}H)^{> k}(x) \rbr ^ {2} \one_{x\in X \andd \tau_{x} \neq 0}\rbs,
\end{split}
\end{align}
with the first equality a consequence of $H(x)=\ee{I}H(x)$ in case $\tau_{x}=0$.

We first bound the second summand appearing on the RHS of~\eqref{eq:eating3}. Using projection arguments (Claim~\ref{clm:projection_shrink}), plus the fact that we may take the $(\cdot)^{>k}$ operator inside, analogously to~\eqref{Aux5.1}, we obtain:
\begin{align}\label{Eq:Aux5.6}
\begin{split}
\be_{x} \lbs \lbr H^{> k}(x) - (\ee{I}H)^{> k}(x) \rbr ^ {2} \rbs
& = \be_{x} \lbs \lbr \lbr \lbr \one - \ee{I} \rbr \circ D_{P} \rbr \lbr h^{>k}\rbr\rbr ^{>k} (x)^{2} \rbs
\\
& \leq \normts{h^{>k}} \leq \eps.
\end{split}
\end{align}
Now we bound the first summand appearing on the RHS of~\eqref{eq:eating3}. Using the same $(a+b)^2 \leq 2(a^2+b^2)$ argument as before, we obtain, for any $x\in X$,
\begin{align*}
\tau_x &\leq \var \lbr H_{\bar{I} = x} \rbr \\
       &\leq 2 \underbrace{\be_{y} \lbs \given{\lbr H^{\leq k}(y) - (\ee{I}H)^{\leq k}(y) \rbr ^ {2}}{y\equiv_{\bar{I}}x} \rbs}_{\gamma_{x}} +
	2 \underbrace{\be_{y} \lbs \given{\lbr H^{> k}(y) - (\ee{I}H)^{> k}(y) \rbr ^ {2}}{y\equiv_{\bar{I}}x} \rbs}_{\delta_{x}}.
\end{align*}
In particular, if $x \in X$, then either $\gamma_{x} > \tau_{x}/4$ or else $\delta_{x} > \gamma_{x}$. Either way,
\begin{equation}\label{Eq:Aux5.7}
\gamma_{x} \leq \frac{4}{\tau_{x}}\gamma_{x}^2 + \delta_{x}.
\end{equation}
Hence, we can bound the first summand in the RHS of~\eqref{eq:eating3} as follows:
\begin{align}\label{Eq:Aux5.8}
\begin{split}
	&\be\left[(H^{\leq k}-\ee{I}H^{\leq k})(x)^{2} \one_{x\in X \wedge \tau_{x} \neq 0}\right]
	\\
	& \qquad\qquad\qquad\qquad\qquad \stackrel{(a)}{\leq}
	4\be\left[\frac{(H^{\leq k}-\ee{I}H^{\leq k})(x)^{4}}{\tau_x}\one_{x\in X \wedge \tau_{x} \neq 0}\right] +
	\be\left[(H^{> k}-\ee{I}H^{> k})^{2}\right]
	\\
	& \qquad\qquad\qquad\qquad\qquad \stackrel{(b)}{\leq}
	4\sqrt{\be \lbs (1/\tau_x^{2}) \cdot \one_{\tau_{x} \neq 0}\rbs} \sqrt{\be \lbs (H^{\leq k}-\ee{I}H^{\leq k})^{8} \rbs} + \eps
	\\
	& \qquad\qquad\qquad\qquad\qquad \stackrel{(c)}{\leq}
	\sqrt{\be \lbs (1/\tau_x^{2}) \cdot \one_{\tau_{x} \neq 0}\rbs} \sqrt{p ^ {-O(k)} \be \lbs (H^{\leq k}-(\ee{I}H)^{\leq k})^{2} \rbs^{4}} + \eps.
\end{split}
\end{align}
\noindent Inequality $(a)$ follows from~\eqref{Eq:Aux5.7}, and the inequality $\be[Z]^2\leq \be[Z^2]$, applied as
\[
	\be_{x}\lbs \frac{\one_{x\in X \wedge \tau_{x} \neq 0}}{\tau_{x}}\gamma_{x}^{2} \rbs \leq \be_{x} \lbs \frac{\one_{x\in X \wedge \tau_{x} \neq 0}}{\tau_{x}} \be_{y}\lbs \given{\lbr H^{\leq k} - (\ee{I}H)^{\leq k} \rbr(y) ^ {4}}{y \equiv_{\bar{I}} x}\rbs \rbs.
\]
Inequality $(b)$ follows from the Cauchy-Schwarz inequality $\be[A\cdot B] \leq \sqrt{\be[A^{2}]} \cdot \sqrt{\be[B^{2}]}$ and~\eqref{Eq:Aux5.6}, and Inequality $(c)$ is an application of Lemma~\ref{lem:bonami_beckner_slice}.

\medskip Before we proceed to bound the RHS of~\eqref{Eq:Aux5.8}, we note that we can assume $2|I|>K$, where $K=O(k^2/p)$, as otherwise, we can deduce~\eqref{eq:eating_var_bound} easily even without the above argument. Indeed, Lemma~\ref{lem:var_bound_by_larger} applied on $H$ gives
\begin{equation}\label{Eq:aux510}
	\forall i \in I: \qquad \vv{I}(H)\leq \vv{I\sm \{i\}}(H)+2\eps,
\end{equation}
which inductively implies $\vv{I}(H)\leq 2|I|\eps$, as required in the case $2|I| \leq K$. (Notice the involved derivatives of $H$ are indeed bounded by $2\eps$, as we assumed every $(l+1)$-st order derivative of $h$, containing $P$, is small, and $I\cap \bigcup P = \es$.) Hence, from now on we assume $2|I|>K$.

\medskip In order to bound the RHS of~\eqref{Eq:Aux5.8}, we first bound the term $\be \lbs (1/\tau_x^{2})\cdot \one_{\tau_{x} \neq 0}\rbs$, as follows.
\begin{align*}
\begin{split}
	\be[(1/\tau_{x}^{2}) \cdot \one_{\tau_{x} \neq 0}]
	& %\stackrel{(d)}{=}
	= 25 \sum_{s=1}^{|I|-1} \pr\lbs \sum_{i \in I} x_{i}=s\rbs \nu(s/|I|)^{-c_{2}k}
	\\
	& \leq
	50 \sum_{s=1}^{|I|/2} \pr\lbs \sum_{i \in I} x_{i}=s\rbs (|I|/s)^{c_{2}k}
	\\
	& \leq
	50\pr\lbs \sum_{i \in I} x_{i}\geq p|I|/2\rbs (2/p)^{c_{2}k} + 50 \pr\lbs \sum_{i \in I} x_{i}\leq p|I|/2\rbs |I|^{c_{2}k}
	\\
	& \stackrel{(d)}{\leq}
	50 (2/p)^{c_{2}k} + 50 \exp(-p|I|/8) \cdot |I|^{c_{2}k}
	\\
	& \leq
	50 (2/p)^{c_{2}k} + (1/p)^{c_{2}k},
\end{split}
\end{align*}
Inequality $(d)$ uses Lemma~\ref{lem:neg_chernoff}, and the last inequality uses the assumption that $|I|>K/2$ with $K=O(k^2/p)$. We thus obtained
\begin{equation}\label{Eq:Aux5.9}
\be\lbs 1/\tau_{x}^{2} \cdot \one_{\tau_{x} \neq 0} \rbs \leq p^{-O(k)}.
\end{equation}
To conclude the bounding of the RHS of~\eqref{Eq:Aux5.8}, we bound the term $\be \lbs (H^{\leq k}-(\ee{I}H)^{\leq k})^{2} \rbs^{4}$. As we assumed $2|I|>K$, inequality~\eqref{Eq:aux510}, together with the inductive assumption $\vv{I \sm \{i\}}(H) \leq K\eps$, implies $\vv{I}(H) \leq (K+2)\eps$. Using this, along with a projection argument (Claim~\ref{clm:projection_shrink} above), we obtain
\begin{equation}\label{Eq:Aux5.10}
\be \lbs (H^{\leq k}-(\ee{I}H)^{\leq k})^{2} \rbs^{4} \leq ((K+2)\epsilon)^4.
\end{equation}
Combining~\eqref{Eq:Aux5.8},~\eqref{Eq:Aux5.9}, and~\eqref{Eq:Aux5.10}, we get
\begin{align}\label{Eq:Aux5.11}
\begin{split}
\be\left[(H^{\leq k}-\ee{I}H^{\leq k})^{2}\one_{x\in X}\right] &\leq \sqrt{\be \lbs 1/\tau_x^{2}\rbs} \sqrt{p ^ {-O(k)} \be \lbs (H^{\leq k}-(\ee{I}H)^{\leq k})^{2} \rbs^{4}} + \eps \\
&\leq p^{-O(k)}(K+2)^{2}\eps^{2} + \eps \leq 2\eps,
\end{split}
\end{align}
where the last inequality follows from our choice of $K\ll p^{-O(k)}$, and our assumption $\eps < p^{c_{1}k}$.

\medskip Combining Inequalities~\eqref{eq:eating3},~\eqref{Eq:Aux5.6} and~\eqref{Eq:Aux5.11}, we obtain
\begin{equation}\label{Eq:Aux5.12}
\be_{x} \lbs \lbr H(x) - \ee{I}H(x) \rbr ^{2} \one_{x\in X}\rbs \leq 2(\eps + 2\eps) = 6\eps,
\end{equation}
thus completing the second sub-step.

\medskip Finally, combining the above estimates for $x \notin X$ (i.e., Inequality~\eqref{Eq:Aux5.5} above) and for $x \in X$ (i.e., Inequality~\eqref{Eq:Aux5.12} above), and plugging into Equation~\eqref{eq:eating_bound0}, we deduce $\vv{I}(H) \leq 8\eps$, which clearly implies the desired Inequality~\eqref{eq:eating_var_bound} inductively, also in the case $2|I| > K$, provided we choose $K\geq 8$.

\medskip \noindent \textbf{Step~2: Constructing the approximating function.} In Step~1 we proved that for any $l$-tuple $P$, we have $\vv{J}(D_{P}h) \leq K\eps$, where $J=[n] \sm \bigcup P$.

Let
\[
	a_{P}= \frac{\be \lbs (D_{P}h) \cdot \Psi_{P}\rbs}{\Pr_{x}[\Psi_{P}(x)\neq 0]} \in\reals.
\]
We claim that
\[
\vv{J} (D_{P}h) = \be \lbs (D_{P}h - a_{P} \cdot \Psi_{P})^{2} \rbs.
\]
Clearly, it is enough to prove that $\ee{J}(D_{P}h) = a_{P}\Psi_{P}$. This indeed holds, since

\begin{align*}
\begin{split}
	\ee{J}(D_{P}h)(x)
	   & = \be_{\pi \in S_{J}} \be_{T \subseteq P} \lbs (-1)^{|T|} h((x^{T})^{\pi}) \rbs
	\\ & = \be_{\pi \in S_{J}} \be_{T \subseteq P} \lbs h((x^{T})^{\pi}) \Psi_{P}((x^{T}) ^ {\pi}) \Psi_{P}(x) \rbs
	\\ & = \be_{y:\Psi_{P}(y)\neq 0}\lbs h(y) \cdot \Psi_{P}(y)\rbs \cdot \Psi_{P}(x)
	\\ & = \be_{y:\Psi_{P}(y)\neq 0}\ \be_{T\sub P} \lbs h(y^{T}) \cdot \Psi_{P}(y^{T})\rbs\cdot \Psi_{P}(x)
	\\ & = \be_{y:\Psi_{P}(y)\neq 0} \lbs (D_{P}h)(y)\cdot \Psi_{P}(y) \rbs \cdot \Psi_{P}(x)
	\\ & = a_{P}\cdot \Psi_{P}(x).
\end{split}
\end{align*}

Recall that $P$ is of order $l$ and hence $D_{P}h$ is $2^{-l}\integers$-valued. Similarly to Equation~\eqref{eq:round_not_so_bad}, letting $c_{P} \ddd 2^{-l}\lfloor 2^{l} a_{P} \rceil$, we find
\[
	\forall x\colon (D_{P}h(x)-c_{P}\Psi_{P}(x))^2 \leq 4 (D_{P}h(x)-a_{P}\cdot\Psi_{P}(x))^2,
\]
and consequently, for every $P$ we have $$\be\lbs \lbr D_{P}h-c_{P}\Psi_{P} \rbr^{2} \rbs \leq 4K\eps.$$
Since $c_{P}\in 2^{-l}\integers$, Proposition~\ref{lem:low_level_attain_derivatives} guarantees the existence of a function $\sifunc{g}$ of degree $\leq l$ that satisfies $D_{P}g =c_{P}\Psi_{P}$ for all \textbf{shifted sorted} $l$-tuples $P$. We shall see that $g$ is the function Lemma~\ref{lem:eating_derivatives} seeks for; specifically, we confirm that
\[
\normts{D_{P}(h-g)} \leq 2 \eps
\]
holds for all $l$-tuples $P$.
		
\medskip Indeed, Lemma~\ref{lem:derivative_dichotomy} asserts that for \emph{every} $l$-tuple $P$, we either have $\normts{D_{P}(h-g)}\leq 2\eps$, or else $\normts{D_{P}(h-g)} \geq p^{O(k)}$. The latter case is impossible for \emph{shifted sorted} $l$-tuples, since $4K\eps < p^{O(k)}$, and hence, we must have
\[
\normts{D_{P}(h-g)} \leq 2\eps.
\]
for any shifted sorted $P$. The second part of Proposition~\ref{lem:derivative_generating_set} then implies that actually all $l$-tuples $P$ satisfy $\normts{D_{P}(h-g)} \leq 2 \eps$, concluding the proof of Lemma~\ref{lem:eating_derivatives}.
\end{proof}

\subsection{Deriving the Kindler-Safra theorem from Theorem~\ref{thm:main}}
\label{sec:sub:reduction}
The following proposition shows that Theorem~\ref{thm:main} qualitatively implies the Kindler-Safra theorem~\cite{KS} by a `blackbox' reduction. Comparing Theorems~\ref{thm:sharp_ks} and~\ref{thm:sharp}, one can see that this implication is not the most quantitatively efficient, as the bounds we get in the case of the slice are worse (though, only in the lower order term) than the ones we obtain directly in the discrete-cube case.

For the sake of simplicity, we deduce the Kindler-Safra theorem only for the uniform measure on the discrete cube. In a similar way, the theorem for the biased measure $\mu_p$ can be deduced from Theorem~\ref{thm:main} for the slice $\slice{n}{p}$, for any $0 \ll p \ll 1$. We note that similar arguments were presented, e.g., in~\cite{DFH18,Filmus16b}.
\begin{proposition}\label{Prop:Reduction}
There exists a universal constant $c$, such that the following holds. Let $k \in \pintegers$. Suppose that any $\func{f'}{\binom{[n]}{n/2}}{\zo}$ satisfying $W^{>k}(f') \leq \eps$ is $\alpha(\eps)$-close to a degree-$k$ function $\func{g'}{\binom{[n]}{n/2}}{\zo}$, for some continuous function $\func{\alpha}{[0,1]}{[0,1]}$.

Then, if $\alpha(\eps) < 2^{-ck}$, then any $\func{f}{\zo^{n}}{\zo}$ satisfying $W^{>k}(f) \leq \eps$ is $\alpha(\eps)$-close to a degree-$k$ function $\func{g}{\zo^{n}}{\zo}$.
\end{proposition}
\begin{proof}
	Let $\func{f}{\zo^{n}}{\zo}$ have $W^{>k}(f) \leq \eps$. Choose a large $m \in 2\integers$ and consider the function $\func{f'}{\binom{[m]}{m/2}}{\zo}$ defined by $f'(x) = f\lbr \restrict{x}{[n]}\rbr$. By Definition~\ref{def:level_k}, the functions $\{\AND_{T}\}$ defined over the slice are of degree $\leq |T|$. In particular, a projection argument implies
	\begin{equation}\label{eq:imply_aux1}
		W^{>k}(f')\leq \be_{x\sim \binom{[m]}{m/2}} \lbs {f^{>k}\lbr \restrict{x}{[n]}\rbr ^{2}} \rbs,
	\end{equation}
	where $f^{>k}$, is, as usual, a function defined over $\zo^{n}$. The key observation is that as $m\to \infty$, the distribution of the random variable $y=\restrict{x}{[n]}$ for $x\sim \binom{[m]}{m/2}$ approaches the uniform measure on $\zo^{n}$. This is because $y$ attains every value in $\zo^{n}$ with some probability $q$ satisfying
\[
\binom{m-n}{m/2-n}/\binom{m}{m/2} \leq q \leq \binom{m-n}{(m-n)/2}/\binom{m}{m/2},
\]
which is $2^{-n}\lbr 1 \pm O(n/\sqrt{m}) \rbr$ for $m > n^{2}$.

This, together with~\eqref{eq:imply_aux1}, implies $W^{>k}(f') \leq (1+o_m(1))\eps$, where $o_m(1)$ denotes a quantity that tends to 0 as $m \to \infty$. Hence, $f'$ is $(1+o(1))\alpha(\eps)$-close to a degree-$k$ function $\func{g'}{\binom{[m]}{m/2}}{\zo}$.
	
\medskip An important property of $g'$ is that `it does not depend on coordinates outside $[n]$', that is, $D_{ij}g'\equiv 0$ for every $i,j\in [m]\sm [n]$. To see this, let $i,j \in [m]\sm [n]$, and consider $g'^{(ij)}$. Evidently, by definition of $f'$, one has $f' = f'^{(ij)}$, and therefore,
	\[
		\pr[f' \neq g'^{(ij)}] = \pr[f'^{(ij)} \neq g'] = \pr[f' \neq g'] \leq (1+o(1))\alpha(\eps).
	\]
	Combining this with the triangle inequality, we see that $g', g'^{(ij)}$ are both degree-$k$ functions, which are $3\alpha(\eps)$-close to each other. However, Lemma~\ref{lem:integer_to_boolean} states that either $\pr[g' \neq g'^{(ij)}] > 2^{-ck}$ (for some universal constant $c$) or $g'=g'^{(ij)}$. Since, by assumption, $\alpha(\eps)$ is small enough, we must have $g' = g'^{(ij)}$, or equivalently, $D_{ij}g' \equiv 0$, as asserted.

\medskip Finally, we may define $\func{g}{\zo^{n}}{\zo}$ by
	\[
		g(z)=g'(z_{1}, \ldots, z_{n}, 1-z_{1}, \ldots, 1-z_{n}, 0, 1, 0, 1,\ldots, 0,1),
	\]
	so that $g$ is a degree-$k$ function. Using the invariance of $g'$ to permuting coordinates outside $[n]$, and the fact that the random variables $x\sim \binom{[m]}{m/2}$ and $y=\restrict{x}{[n]}$ are statistically close, one finds that $f$ and $g$ are $(1+o_m(1))\alpha(\eps)$-close, as required. (Taking $m\to \infty$ concludes the proof, as there is only a finite number of Boolean functions over the cube $\zo^{[n]}$.)
\end{proof}

\section{Refinement of the Kindler-Safra Theorem}
\label{sec:sharp}
It turns out that a simple black-box argument using the level-$k$ inequalities can strengthen the Kindler-Safra theorem, and obtain Theorem~\ref{thm:sharp_ks}. Let us recall the level-$k$ inequality for the discrete cube (\cite[Section 9.5]{O'Donnell14}), and then prove Theorem~\ref{thm:sharp_ks}.
\begin{theorem}
Let $\func{h}{\zo^{n}}{\ozo}$ and $k \in \pintegers$ be such that $\mathbb{E}|h| = \alpha \leq \exp(-k/2)$.
Then
\[
	W^{\leq k}(h)\leq \alpha^{2} \lbr \frac{2e}{k}\ln(1/\alpha) \rbr^{k}.
\]
\end{theorem}

\begin{proof}[Proof of Theorem~\ref{thm:sharp_ks}]
Let $\func{g}{\zo^{n}}{\zo}$ be the approximation of $f$ guaranteed by the Kindler-Safra
Theorem (\ref{thm:original_ks}), i.e., $g$ is of degree $\leq k$ and $\pr[f\neq g]\leq O(\eps)$.

Consider $h=f-g$. One has $\func{h}{\zo^{n}}{\ozo}$, and $\be |h|\leq \pr[f\neq g] \leq O(\eps)$. Setting $\alpha=\be|h|$,
we have
\[
\alpha= \be \lbs h^{2} \rbs=W^{\leq k}(h) + W^{>k}(h) \leq \alpha^{2}\lbr \frac{2e}{k}\log(1/\alpha) \rbr^{k}+\eps,
\]
which implies the assertion~\eqref{eq:sharp_ks_assertion} since $\alpha \leq O(\eps)$.
\end{proof}

One can prove Theorem~\ref{thm:sharp} in a similar way, using Theorem~\ref{thm:main} together with the level-$k$ inequalities for the slice (i.e.,  Lemma~\ref{lem:level_k_slice}). Due to the involved less attractive expressions, the proof is omitted.

\medskip We conclude this paper with an example which demonstrates the tightness of Theorem~\ref{thm:sharp_ks}, up to lower-order terms. In order to present the example, we need the following claim.
\begin{claim}\label{clm:biased_halfspace_weight}
	Consider the sequence of Linear Threshold Functions $h_{n}:\zo^{n} \to \spm$, defined by $h_{n}(x) = 1 - 2\cdot \one \{\sum_{i\in [n]} x_{i} > n/2 + t\sqrt{n}/2 \}$ with $t > 0$, so that $\be[h_{n}] \xrightarrow{n\to \infty} 1-\delta > 1/2$. Then, for all $k \in \pintegers$ with $4k+1 < \log(1/\delta)$ and $n$ sufficiently large with respect to $\delta$, we have
	\[
		W^{=k}(h_{n}) \geq \delta^{2} \log(1/\delta)^{k} / O(k)^{k}.
	\]
\end{claim}
\begin{proof}[Proof-Sketch]
The case $k=1$ is trivial, hence we assume $k \geq 2$. It is not hard to see that as $n\to \infty$, we have $1- \be[h] \to \Phi(t)$, where $\Phi(t) = \int_{-\infty}^{t} \phi(x) \mathrm{dx}$ and $\phi(x)=\exp(-x^2/2)/\sqrt{2\pi}$ are the cumulative distribution function (c.d.f.) and probability density function (p.d.f.) of a $\mathcal{N}(0,1)$ random variable. Hence we may assume $\Phi(t) = 1-\delta$. In turn, properties of discrete differentiation imply that $\widehat{h}([k]) \xrightarrow{n\to \infty} \Phi^{(k)}(t) / \sqrt{n}^{k}$. Thus,
	\[
		W^{=k}(h) \approx \frac{\binom{n}{k}}{n^{k}} \Phi^{(k)}(t)^{2} \to \Phi^{(k)}(t)^{2} / k!.
	\]
	Noticing that for $t > 0$,
	\[
		\log(1/\delta) = \log(1/(1-\Phi(t))) < t^{2}+1 \leq (t+1)^{2},
	\]
	we are left with arguing $|\Phi^{(k)}(t)| > \Omega(1+t)^{k} (1-\Phi(t))$ for $t>\sqrt{4k}$.

It is well known that $\Phi^{(k)}(t) = \phi(t) He_{k-1}(t)$, where $He_{k-1}(t)$ is a Hermite polynomial (see~\cite{Szego75}). Moreover, it may be verified that $\phi(t) \geq \frac{3}{4}(1+t)(1-\Phi(t))$. Thus,
\[
|\Phi^{(k)}(t)| \geq \frac{3}{4}(1+t) |He_{k-1}(t)| (1-\Phi(t)).
\]
It only remains to show that $|He_{k-1}(t)| \geq \Omega(1+t)^{k-1}$ whenever $t \geq 2\sqrt{k}$. For this, notice that Hermite polynomials satisfy $He_{i}^{(j)} = \frac{i!}{(i-j)!}He_{i-j}$. Corollary 2.3 of~\cite{ADGR04} implies that the (real) roots of $He_{i}$ are upper bounded by $\sqrt{2i}$. Using the fact that $\lim_{x\to\infty} He_{i}(x)=\infty$, we deduce $\forall j: He_{i}^{(j)}(\sqrt{2i}) \geq 0$. Noting that $He_{i}^{(i)}=i!$, we inductively infer that the polynomial $He_{i}(t)$ has all derivatives larger than the corresponding ones of $(t-\sqrt{2i})^{k}$, for all $t \geq \sqrt{2i}$. In particular, $He_{i}(t) \geq (t-\sqrt{2i})^{i}$ whenever $t > \sqrt{2i}$. We conclude by taking $i=k-1$ and noting that $(t-\sqrt{2k-2})^{k-1} = \Omega(1+t)^{k-1}$ for $t \geq 2\sqrt{k-1}$.
\end{proof}

\begin{example}
Consider the functions $\func{f}{\zo^{n}}{\zo}$ and $\func{g}{\zo^{n}}{\spm}$, defined by
	\begin{equation*}
	\begin{split}
	&g(x) = (-1)^{\sum_{i \leq k} x_{i}}.\\
	&f(x) =
			\begin{cases}
				(1+g(x))/2, &\text{if }\sum_{j=k+1}^{n} x_{j}g(x)^{j} < (n-k+t\sqrt{n-k})/2\\
				(1-g(x))/2, &\text{otherwise.}
			\end{cases}
	\end{split}
	\end{equation*}
	The probability $\delta=\pr_{x}\lbs f(x) \neq (1+g(x))/2 \rbs$ can be arranged to be as small as we wish, by controlling $t\sim \sqrt{2\log(\delta)}$. In particular, in the case $\delta < 2^{-k-1}$, the degree-$k$ function best approximating $f$ is just $(1+g(x))/2$, since the distance between any two degree-$k$ Boolean-valued functions is at least $2^{-k}$ (analogously to Lemma~\ref{lem:integer_to_boolean}). Moreover, inspecting the Fourier expansion of $f$, one finds out
	\[
		\epsilon \ddd W^{>k}(f) = \var(f) - W^{\leq k}(f) \leq \delta - \delta^{2} \log(1/\delta)^{k} / O(k)^{k},
	\]
	as $n\to \infty$. To see this, notice $\var(f)=\delta$, and $f(x)=1/2+g(x) \cdot h(x)$ with
	\[
		h(x)= 1-2 \cdot \one \lbc \sum_{j=k+1}^{n} x_{j}g(x)^{j} > (n-k+t\sqrt{n-k})/2 \rbc,
	\]
	so in order to estimate $W^{\leq k}(f)$, we suffice to understand $W^{\leq k}(g \cdot h)$. We note that $h(x)$ is actually a biased halfspace $h'$ applied on the variables $\{x_{j} g(x)^{j}\}_{j\in[n]\sm [k]}$. Such regular halfspaces with $\be[h']=1-4\delta$ have $W^{\leq k}(h')-W^{0}(h') \geq \delta^{2} \log(1/\delta)^{k} / O(k)^{k}$, as was shown in Claim~\ref{clm:biased_halfspace_weight} (assuming $\delta \leq 2^{-6k-6}$). Each monomial $c x^{S}$ of $h'$ corresponds to either the monomial $c x^{S}$, or to $c x^{S} g(x)$, in $h(x)$. Notice also that about half of such monomials with $0\neq |S|\leq k$ correspond to $c x^{S} g(x)$ (depending on the parity of $\sum_{j \in S} j$). This means
	\[
		W^{\leq k}(g \cdot h) \gtrsim \frac{1}{2} \cdot (W^{\leq k}(h')-W^{0}(h')) \geq \delta^{2} \log(1/\delta)^{k} / O(k)^{k}.
	\]
	We thus obtain
	\[
		\delta > \epsilon + \epsilon^{2} \log(1/\epsilon)^{k} / O(k)^{k},
	\]
	demonstrating that Theorem~\ref{thm:sharp_ks} cannot be improved much. (Recall $O(k)^{k}/k! = 2^{O(k)}$.)
\end{example}

\section*{Acknowledgements}
We are deeply grateful to Yuval Filmus for numerous helpful comments and suggestions, and wish to thank Guy Kindler for explaining to us some aspects of his work~\cite{KS}.


\begin{thebibliography}{99}

\bibitem{AK} R. Ahlswede and L. H. Khachatrian, The complete intersection theorem for systems of finite
sets, {\it Eur. J. Combin.} \textbf{18} (1997), pp. 125--136.

\bibitem{ADFS} N. Alon, I. Dinur, E. Friedgut, and B. Sudakov, Graph Products, Fourier Analysis and Spectral Techniques, {\it Geom. Func. Anal.} \textbf{14(5)} (2004), pp.~913--940.

\bibitem{ADGR04} I. Area, D. K. Dimitrov, E. Godoy, and A. Ronveaux, Zeros of Gegenbauer and Hermite polynomials and connection coefficients,
{\it Math. Comp.}, \textbf{73} (2004), pp.~1937--1951.

\bibitem{BNR16} B. Bollob\'{a}s, B. P. Narayanan, and A. M. Raigorodskii, On the stability of the Erd\H{o}s-Ko-Rado theorem, {\it J. Combin. Theory Ser. A}
\textbf{137} (2016), pp. 64--78.

\bibitem{Bonami} A. Bonami, Etude des coefficients Fourier des fonctiones de $L^p (G)$, {\it Ann. Inst. Fourier} \textbf{20} (1970), pp.~335--402.

\bibitem{BKKKL} J. Bourgain, J. Kahn, G. Kalai, Y. Katznelson, and N. Linial, The influence of variables in product spaces, {\it
Israel J. Math} \textbf{77} (1992), no. 1--2, pp.~55--64.

\bibitem{CHS18} J. Chiarelli, P. Hatami, and M. Saks, An asymptotically tight bound on the number of relevant variables in a bounded degree Boolean function, 2018. Available at arxiv:1801.08564.

\bibitem{DT16} S. Das and T. Tran, Removal and stability for Erd\H{o}s--Ko--Rado, {\it SIAM J. Disc. Math.}, \textbf{30(2)} (2016), pp. 1102--1114.

\bibitem{Dinur07} I. Dinur, The PCP theorem by gap amplification, {\it Journal of the ACM}, \textbf{54(3)} (2007), pp.~1--44.

\bibitem{DFH19} I. Dinur, Y. Filmus, and P. Harsha, Analyzing Boolean functions on the biased hypercube via higher-dimensional agreement tests (Extended abstract), proceedings of ACM-SIAM Symposium on Discrete Algorithms (SODA) 2019, pp.~2124--2133.

\bibitem{DFH18} I. Dinur, Y. Filmus, and P. Harsha, Low degree almost Boolean functions are sparse juntas, preprint, 2017. Available at arXiv:1711.09428.

\bibitem{EFP11} D. Ellis, E. Friedgut, and H. Pilpel, Intersecting families of permutations, {\it J. Amer. Math. Soc.} \textbf{24} (2011), pp.~649--682.

\bibitem{EKL16} D. Ellis, N. Keller and N. Lifshitz, Stability versions of Erd\H{o}s-Ko-Rado type theorems, via isoperimetry, {\it J. Eur. Math. Soc.}, to appear. arXiv:1604.02160.

\bibitem{FF14} D. Falik and E. Friedgut, Between Arrow and Gibbard-Satterthwaite: A representation theoretic approach, {\it Israel J. Math} \textbf{201(1)} (2014), pp.~247--297.

\bibitem{Filmus16} Y. Filmus, An orthogonal basis for functions over a slice of the Boolean cube, {\it Electron. J. Combin.} \textbf{23(1)}:P1.23, 2016.

\bibitem{Filmus16b} Y. Filmus, Friedgut-Kalai-Naor theorem for slices of the Boolean cube, {\it Chicago J. Theor. Comput. Sci.} 2016:14, 2016.

\bibitem{FI18} Y. Filmus and F. Ihringer, Boolean constant degree functions on the slice are juntas, 2018. arXiv:1801.06338.

\bibitem{FKMW18} Y. Filmus, G. Kindler, E. Mossel, and K.Wimmer, Invariance principle on the slice, {\it ACM Trans. Comput. Th.} \textbf{10(3)}:11, 2018.

\bibitem{FM16} Y. Filmus and E. Mossel, Harmonicity and invariance on slices of the Boolean cube, proceedings of CCC'2016, pp.~16:1--16:13, 2016.

\bibitem{FM16*} Y. Filmus and E. Mossel, Harmonicity and invariance on slices of the Boolean cube (full version). arXiv:1507.02713, 2015.

\bibitem{Frankl87} P. Frankl, The shifting technique in extremal set theory, in: {\it Surveys in Combinatorics,
Lond. Math. Soc. Lect. Note Ser.} \textbf{123} (1987), pp. 81--110.

\bibitem{FG89} P. Frankl and R. L. Graham, Old and new proofs of the Erd\H{o}s-Ko-Rado theorem. {\it J. Sichuan Univ. Nat. Sci. Ed.}, \textbf{26} (1989), pp.~:112--122.

\bibitem{FT16} P. Frankl and N. Tokushige, Invitation to intersection problems for finite sets, {\it J. Combin. Th.
Ser. A} \textbf{144} (2016), pp.~157--211.

\bibitem{Friedgut98} E. Friedgut, Boolean functions with low average sensitivity depend on few coordinates, {\it Combinatorica}
\textbf{18(1)} (1998), pp.~27--35.

\bibitem{Friedgut08} E. Friedgut, On the measure of intersecting families, uniqueness and
stability, {\it Combinatorica} \textbf{28} (2008), pp. 503--528.

\bibitem{FKN} E. Friedgut, G. Kalai, and A. Naor, Boolean functions whose Fourier transform
is concentrated on the first two levels, {\it Adv. Appl. Math.}, \textbf{29(3)} (2002), pp.~427--437.

\bibitem{GH08} M. Ghandehari and H. Hatami, Fourier analysis and large independent sets in powers of complete graphs, {\it J. Combin. Th. Ser. B} \textbf{98(1)} (2008), pp.~164--172.

\bibitem{GG06} B. Graham and G. R. Grimmett, Influence and sharp-threshold theorems for monotonic measures, {\it Ann. Probab.} \textbf{34} (2006), pp.~1726--1745.

\bibitem{JOW15} J. Jendrej, K. Oleszkiewicz, and J. O. Wojtaszczyk. On some extensions of the FKN theorem, {\it Theory of Comput.} \textbf{11} (2015), pp.~445--469.

\bibitem{KKL} J. Kahn, G. Kalai, and N. Linial, The influence of variables on Boolean functions, proccedings of FOCS'1988, pp.~68--80, 1988.

\bibitem{Kalai-Choice} G. Kalai, A Fourier-theoretic perspective on the Condorcet paradox and Arrow's theorem, {\it Adv. Appl. Math.} \textbf{29(3)} (2002), pp.~412--426.

\bibitem{Karpas17} I. Karpas, Two results on union-closed families, available at arXiv:1708:01434, 2017.

\bibitem{KS} G. Kindler and S. Safra, Noise-resistant Boolean functions are juntas, Manuscript, 2002.

\bibitem{LY98} T.-Y. Lee and H.-T. Yau, Logarithmic Sobolev inequality for some models of random walks, {\it Ann. Probab.} \textbf{26(4)} (1998), pp.~1855--1873.

\bibitem{MO10} A. Montanaro and T. Osborne, Quantum boolean functions, {\it Chicago J. Theor. Comput. Sci.}, 2010.

\bibitem{Maj-Stablest} E. Mossel, R. O'Donnell and K. Oleszkiewicz, Noise stability of functions with low influences:
Invariance and optimality, {\it Annals of Math.} \textbf{175(3)} (2012), pp.~1283--1327.

\bibitem{Nayar14} P. Nayar, FKN theorem on the biased cube, {\it Colloq. Math.} \textbf{137(2)} (2014), pp.~253--261.
2014.

\bibitem{NS94} N. Nisan and M. Szegedy, On the degree of Boolean functions as real polynomials, {\it Comp. Comp.} \textbf{4(4)} (1994), pp.~301--313.

\bibitem{O'Donnell14} R. O'Donnell, Analysis of Boolean functions, Cambridge University Press, 2014.

\bibitem{OW13} 	R. O'Donnell and K. Wimmer, KKL, Kruskal-Katona, and monotone nets, {\it SIAM J. Comput.} \textbf{42(6)} (2013), pp.~2375--2399.

\bibitem{PS97} A. Panconesi and A. Srinivasan, Randomized distributed edge coloring via an extension of the Chernoff-Hoeffding bounds, {\it SIAM J. Comput.} \textbf{26} (1997), pp.~350--368.

\bibitem{Rub12} A. Rubinstein, Boolean functions whose Fourier transform is concentrated on pair-wise disjoint subsets of the inputs, M. Sc. thesis, Tel Aviv University, 2012.

\bibitem{Samorodnitsky16} A. Samorodnitsky, On the entropy of a noisy function, {\it IEEE Trans. Information Theory} \textbf{62(10)} (2016), pp.~5446--5464.

\bibitem{Srinivasan11} M. K. Srinivasan, Symmetric chains, Gelfand-Tsetlin chains, and the Terwilliger algebra of the
binary Hamming scheme, {\it J. Algebr. Comb.} \textbf{34(2)} (2011), pp.~301--322.

\bibitem{Szego75} G. Szeg\"{o}, Orthogonal Polynomials (4th ed.), Amer. Math. Soc. Coll. Publ., vol. 23, Providence, RI, 1975.

\bibitem{Wimmer14} K. Wimmer, Low influence functions over slices of the Boolean hypercube depend on few coordinates, proceedings of CCC'2014, pp.~120--131.
\end{thebibliography}
\end{document}